\documentclass[11pt]{amsart} 

\usepackage{amsmath, amscd, amsthm, amssymb, mathrsfs} 

\DeclareMathOperator{\im}{im}

\DeclareMathOperator{\End}{End}

\DeclareMathOperator{\tr}{tr}

\DeclareMathOperator{\Ric}{Ric}
\DeclareMathOperator{\vol}{vol}

\DeclareMathOperator{\Scal}{Scal}

\newcommand{\R}{\mathbb R}
\newcommand{\C}{\mathbb C}

\newcommand{\diff}{\text{\rm d}}
\newcommand{\del}{\partial}
\newcommand{\delb}{\bar{\del}}

\renewcommand{\P}{\mathbb P}

\theoremstyle{plain}
	\newtheorem{theorem}{Theorem}
	
	\newtheorem{lemma}[theorem]{Lemma}
	\newtheorem{corollary}[theorem]{Corollary}

\theoremstyle{definition}
	\newtheorem{definition}[theorem]{Definition}
	\newtheorem{remark}[theorem]{Remark}
	\newtheorem{remarks}[theorem]{Remarks}
	
	\newtheorem{example}[theorem]{Example}

\theoremstyle{plain}
	\newtheorem*{theorem*}{Theorem}
	\newtheorem*{proposition*}{Proposition}
	\newtheorem*{lemma*}{Lemma}
	\newtheorem*{corollary*}{Corollary}
	\newtheorem*{conjecture*}{Conjecture}
\theoremstyle{definition}
	\newtheorem*{definition*}{Definition}
	\newtheorem*{remark*}{Remark}
	\newtheorem*{remarks*}{Remarks}

\numberwithin{equation}{section}
\numberwithin{theorem}{section}

\begin{document}

\title[Constant scalar curvature metrics]{Constant scalar curvature K\"ahler metrics on fibred complex surfaces}
\author{Joel Fine}
\address{Imperial College\\ Mathematics Department\\ London SW7 2AZ\\ United Kingdom}
\email{\tt joel.fine@imperial.ac.uk}
\date{}

\begin{abstract}
This article finds constant scalar curvature K\"ahler metrics on certain compact complex surfaces. The surfaces considered are those admitting a holomorphic submersion to curve, with fibres of genus at least $2$. The proof is via an adiabatic limit. An approximate solution is constructed out of the hyperbolic metrics on the fibres and a large multiple of a certain metric on the base. A parameter dependent inverse function theorem is then used to perturb the approximate solution to a genuine solution in the same cohomology class. The arguments also apply to certain higher dimensional fibred K\"ahler manifolds.
\end{abstract} 

\maketitle

\section{Introduction}\label{introduction}

This article proves the existence of constant scalar curvature K\"ahler metrics on certain complex surfaces. Let $\pi \colon X \to \Sigma$ be a holomorphic submersion from a compact connected complex surface to a curve, with fibres of genus at least $2$. Topologically, $X$ is a locally trivial surface bundle; analytically, however, the fibrewise complex structure may vary. Examples of such surfaces appear in, for example, \cite{atiyah:tsofb, kodaira:actoias}. 

$X$ is polarised by $\kappa_{r} = -c_{1}(V)-r\pi^{*}c_{1}(\Sigma)$, where $V$ is the vertical tangent bundle and the integer $r$ is sufficiently large. More generally, this is a K\"ahler class when $r$ is non-integral.

The precise result proved here is:
\begin{theorem}\label{main theorem}
If $X$ is a compact connected complex surface admitting a holomorphic submersion onto a complex curve with fibres of genus at least two, then, for all large $r$, the K\"ahler class $\kappa_r$ contains a constant scalar curvature K\"ahler metric.
\end{theorem}

In fact, the arguments apply to certain higher dimensional fibred K\"ahler manifolds, although the conditions are more awkward to state; see Theorem \ref{higher dimensional theorem}.

\subsection{Stability of polarised varieties}

There are many reasons to be interested in constant scalar curvature K\"ahler (cscK) metrics. Perhaps the most important comes from a well known conjecture relating the existence of a cscK metric to the stability of the underlying polarised variety. This is analogous to the Hitchin--Kobayashi correspondence and is due, in both chronological order and degree of generality, to Yau \cite{yau:opig}, Tian \cite{tian:kemwpsc} and Donaldson \cite{donaldson:rogtcga4mt}.

The surfaces covered by Theorem \ref{main theorem} are each a family of stable curves parameterised by a stable base, and so stand a good chance of being stable, at least for a polarisation in which the base is large. This is certainly true in the analogous situation for bundles. 

The algebro-geometric view point is not explored here, except for a few comments relating analytic observations to stability conditions.

\subsection{Projectivised vector bundles}

A result analogous to Theorem \ref{main theorem}, but concerning projectivised vector bundles, is proved by Hong in \cite{hong:rmwchsc} and \cite{hong:chsceorm}. The techniques used there are similar to those applied here. There are two important differences between the situations, however.

In both cases the fibres have cscK metrics. However, in Theorem \ref{main theorem} the fibres have  no nonzero holomorphic vector fields, whilst the fibres of a projectivised vector bundle have holomorphic vector fields which are not Killing. If $(F, \omega)$ is a compact K\"ahler manifold with cscK metric, and $\xi$ is a holomorphic vector field on $F$ which is not Killing, then flowing $\omega$ along $\xi$ gives a non-trivial family of cohomologous cscK metrics. 

From the point of view of the analysis, this situation leads to problems, as is explained later. It is avoided in \cite{hong:rmwchsc} and \cite{hong:chsceorm} by the additional assumption that, whilst there are holomorphic vector fields on the fibres, they are not induced by a global holomorphic vector field. From the point of view of the algebraic geometry, this is related to the stability of the vector bundle which is projectivised.

The other important difference between the two situations is that the fibres of a ruled manifold are rigid, whilst the fibres in Theorem \ref{main theorem} have moduli. This leads to extra considerations in the proof, which are not required in Hong's work.

\subsection{Base genus $0$ or $1$}

Let $X$ satisfy the hypotheses of Theorem \ref{main theorem}. When $\Sigma$ has genus $0$ or $1$, a direct proof of the theorem can be given. 

\begin{theorem}
Let $X$ be a compact connected complex surface admitting a holomorphic submersion $\pi\colon X \to \Sigma$ to a surface of genus $0$ or $1$, for which all the fibres are of genus at least $2$. Then $X$ admits cscK metrics.
\end{theorem}

\begin{proof}
Mapping each fibre of $\pi$ to its Jacobian determines a map $j \colon  \Sigma \to \mathscr A$, where $\mathscr A$ denotes the moduli space of principally polarised abelian varieties (of dimension equal to the genus of the fibres of $\pi$). The universal cover of $\mathscr A$ is the Siegel upper half space $\mathscr S$ which can be realised as a bounded domain in $\C^N$. 

The map $j$ lifts to a holomorphic map from the universal cover of $\Sigma$ to $\mathscr S$. Since $g(\Sigma) = 0$ or $1$, this map must be constant: when $g(\Sigma) = 0$ the map is a holomorphic map from a compact manifold; when $g(\Sigma)=1$ the map is a bounded holomorphic map from $\C$. Hence $j$ is constant. By Torelli's theorem, all the fibres of $\pi$ are biholomorphic.

As the model fibre $S$ of $X \to \Sigma$ has genus at least $2$, its group of biholomorphisms $\Gamma$ is finite. Define a principal $\Gamma$ bundle $P \to \Sigma$ by setting the fibre over $\sigma$ to be the group of biholomorphisms from $\pi^{-1}(\sigma)$ to $S$. Since $P$ is a cover of $\Sigma$, it arises from some representation $\pi_1(\Sigma) \to \Gamma$. Using this representation, 
$$
X = P \times_{\pi_1(\Sigma)} S.
$$

In the case $g(\Sigma)=0$, this gives $X = S \times \P^1$ which admits product cscK metrics. If $g(\Sigma) = 1$, $X$ is a quotient of $S \times \C$ by $\pi_1(\Sigma)$. The product $S \times \C$ admits a natural cscK metric, and with respect to this metric $\pi_1(\Sigma)$ acts by isometries. Hence the metric descends to a metric on $X$ which also has constant scalar curvature.
\end{proof}

From now on, in the proof of Theorem \ref{main theorem} only the case $g(\Sigma) \geq 2$ is considered.

\subsection{Outline of proof}

This section gives an outline of the proof of Theorem \ref{main theorem}. The main technique used is an \emph{adiabatic limit}. This involves studying a family of metrics on $X$ for which the base becomes increasingly large.

The first step is to construct a family of approximate solutions. This is done in Section \ref{approximate solutions}. The motivating idea is that in an adiabatic limit the local geometry is dominated by that of the fibre. The approximate solutions are constructed by fitting together the cscK metrics on the fibres of $X$ and a large multiple of a metric on the base. If the base is scaled by a factor $r$, the metrics on $X$ have scalar curvature which is $O(r^{-1})$ from being constant. Section \ref{approximate solutions} also discusses adjusting these metrics to decrease the error to $O(r^{-n})$ for any positive integer $n$.

The remaining sections carry out the analysis necessary to show that a genuine solution lies nearby. Doing this involves solving a parameter dependent inverse function theorem. As is explained in Section \ref{set-up ift}, such arguments hinge on certain estimates. Two sorts of estimate are used. The first kind involves local analysis; the various constants appearing in Sobolev inequalities are shown to be independent of $r$ (see Section \ref{local analysis}). The second kind of estimate involves global analysis, specifically a lower bound for the first eigenvalue of the linearisation of the scalar curvature equation as is found in Section \ref{inverse}.

Before the start of the proof itself, Section \ref{general properties} collates some background information concerning straightforward properties of scalar curvature and its dependence on the K\"ahler structure. In particular it states precisely how various geometric objects are uniformly continuous with respect to the metric used to define them. These results are standard and widely known, although not stated in the literature in quite the form used here.

\subsubsection*{Acknowledgements}

I would like to thank my thesis supervisor Simon Donaldson for suggesting this problem and for sharing his insights with me. I would also like to thank Dominic Joyce for clarifications and corrections to the original text.

\section{General properties of scalar curvature}
\label{general properties}

\subsection{The scalar curvature map}

K\"ahler metrics on a compact K\"ahler manifold $X$ and in a fixed cohomology class are parametrised by K\"ahler potentials: any other K\"ahler form cohomologous to a given one $\omega$ is of the form
$\omega_\phi = \omega + i \delb \del \phi$ for some real valued $\phi \in C^\infty$. 

The scalar curvature of a K\"ahler metric can be expressed as the trace of the Ricci form with respect to the K\"ahler form: $\Scal(\omega)\, \omega^n = n \rho \wedge \omega^{n-1}$,
where $n$ is the complex dimension of $X$. The equation studied in this article is
$$
\Scal(\omega_\phi) = \text{const.}
$$
This is a fourth order, fully nonlinear, partial differential equation for $\phi$. 

Scalar curvature determines a map $S\colon\phi \mapsto \Scal(\omega_\phi)$ defined, initially, on an open set $U \subset C^{\infty}$.

\begin{lemma}
Let $V$ denote the $L^p_{k+4}$-completion of $U$. $S$ extends to a smooth map $S : V \to L^p_k$ whenever $(k+2)p - 2n > 0$. Its derivative at $0$ is given by 
\begin{equation}\label{derivative of scalar curvature}
L(\phi) 	= 	\left(\Delta^2 - S(0) \Delta\right)\phi
				+ n(n-1)\frac{i \delb\del \phi \wedge \rho \wedge \omega^{n-2}}
								{\omega^n}.
\end{equation}
\end{lemma}

\begin{proof}
Scalar curvature is analytic in the metric, and hence $S$ extends to a smooth map from $L^p_{k+4}$ provided that $i \delb\del\phi$ is continuous, \emph{i.e.}\ for $\phi \in L^{2}_{k+4}$ where $(k+2)p -2n>0$.

To compute its derivative, let $\omega_t = \omega + t i \delb\del \phi$. The corresponding metric on the tangent bundle is $g_t = g +t \Phi$ where $\Phi$ is the real symmetric tensor corresponding to $i\delb\del\phi$. The Ricci form is given locally by
$$
\rho_t = \rho + i \delb \del \log \det \left( 1 + tg^{-1}\Phi \right).
$$
Now $\tr(g^{-1} \Phi) = \Lambda (i\delb\del\phi) = \Delta \phi$. Hence, at $t=0$,
$$
\frac{\diff \omega}{\diff t} = i\delb\del\phi, 
\quad \quad
\frac{\diff \rho}{\diff t} = i\delb\del (\Delta \phi).
$$
The result follows from differentiating the equation $S\omega^n = n \rho\wedge\omega^{n-1}$.
\end{proof}

\begin{remarks}\label{remarks on L}
In flat space $L = \Delta^2$; in general, the leading order term of $L$ is $\Delta^2$. It follows immediately that $L$ is elliptic, with index zero.

It follows, either from the formula for $L$, or from the fact that $\int S(\phi) \omega_\phi^n$ is constant, that
$$
\int L(\phi)\, \omega^n
=
-\int  \phi \Delta S(0) \, \omega^n.
$$
In particular, for cscK metrics, $\im L$ is $L^2$-orthogonal to the constant functions.
\end{remarks}

\begin{lemma}\label{regularity}
Let $k \geq 2$. If $\phi \in C^{k,\alpha}$ satisfies $S(\phi)\in C^{k,\alpha}$ then $\phi \in C^{k+4,\alpha}$.
\end{lemma}

\begin{proof}
$S(\phi) = \Delta_{\phi}V$ where $\Delta_{\phi}$ is the $\omega_{\phi}$-Laplacian and $V = -\log\det(g + \Phi)$ with $\Phi$ being the real symmetric tensor associated to the $(1,1)$-form $i \delb\del \phi$.

Since $\phi \in C^{k, \alpha}$, $\Delta_{\phi}$ is a linear second order elliptic operator with coefficients in $C^{k-2, \alpha}$. By elliptic regularity (see \emph{e.g.}\ \cite{aubin:naommae}, page 87) and assumption on $S(\phi)$, $V \in C^{k, \alpha}$.

The map $\phi \mapsto -\log\det(g + \Phi)$ is nonlinear, second order and elliptic. Such maps also satisfy a regularity result (see \emph{e.g.}\ \cite{aubin:naommae}, page 86), hence $\phi \in C^{k+2, \alpha}$. 

This implies that $\Delta_{\phi}$ has $C^{k,\alpha}$ coefficients meaning that in fact $V \in C^{k+2, \alpha}$ and so $\phi \in C^{k+4,\alpha}$.
\end{proof}

\begin{example}[High genus curves]
\label{L for high genus curves}
For the hyperbolic metric on a high genus curve, $S(0) =-1$. Hence the above lemma gives $L = \Delta^2 + \Delta$. The kernel of $L$ is precisely the constant functions. As $L$ is self-adjoint, it is an isomorphism when considered as a map between spaces of functions with mean value zero.
\end{example}

In the above discussion of scalar curvature, the underlying complex manifold $(X, J)$ is regarded as fixed, whilst the K\"ahler form $\omega$ is varying. An alternative point of view is described in \cite{donaldson:rogtcga4mt}. There the symplectic manifold $(X, \omega)$ is fixed, whilst the complex structure varies (through complex structures compatible with $\omega$). The two points of view are related as follows.

Calculation shows that, on a K\"ahler manifold $(X, J, \omega)$, 
$$
i\delb\del \phi = \mathscr L_{\nabla \phi} \omega.
$$
Hence, infinitesimally, the change in $\omega$ due to the K\"ahler potential $\phi$ is precisely that caused by flowing $\omega$ along $\nabla \phi$. Flowing the K\"ahler structure back along $-\nabla \phi$ restores the original symplectic form, but changes the complex structure by $-\mathscr L _{\nabla \phi} J$. This means that the two points of view (varying $\omega$ versus varying $J$) are related by the diffeomorphism generated by $\nabla \phi$.

The following result, on the first order variation of scalar curvature under changes in complex structure, is proved in \cite{donaldson:rogtcga4mt}. The operator 
$$
\mathscr D :C^\infty(X) \to \Omega^{0,1}(TX)
$$
is defined by $\mathscr D \phi = \delb \nabla \phi$, where $\delb$ is the $\delb$-operator of the holomorphic tangent bundle. The operator $\mathscr D^*$ is the formal adjoint of $\mathscr D$ with respect to the $L^2$-inner product determined by the K\"ahler metric.

\begin{lemma}
An infinitesimal change of $-\mathscr L_{\nabla \phi}J$ in the complex structure $J$ causes an infinitesimal change of $\mathscr D^*\mathscr D \phi$ in the scalar curvature of $(X, J, \omega)$.
\end{lemma}

Taking into account the diffeomorphism required to relate this point of view to that in which $\omega$ varies gives the following formula for the linearisation of scalar curvature with respect to K\"ahler potentials:

\begin{lemma}
\begin{equation}\label{derivative of scalar curvature 2}
L(\phi) = \mathscr D^*\mathscr D \phi + \nabla \Scal \cdot \nabla \phi
\end{equation}
\end{lemma}

Alternatively this can be deduced from equation (\ref{derivative of scalar curvature}) and a Weitzenb\"ock-type formula relating $\mathscr D^{*}\mathscr D$ and $\Delta^{2}$.

If the scalar curvature is constant, then $L = \mathscr D^* \mathscr D$. In particular $\ker L = \ker \mathscr D$ consists of functions with holomorphic gradient. If $X$ has constant scalar curvature and no holomorphic vector fields, then $\ker L = \R$. Since $L$ is also self-adjoint it has index zero and so is an isomorphism between spaces of functions with mean value zero. This generalises Example \ref{L for high genus curves}.

\subsection{Dependence on the K\"ahler structure}

This section proves that the scalar curvature map is uniformly continuous under changes of the K\"ahler structure. The arguments are straightforward, this section simply serves to give a precise statement of estimates that will be used later. 

\subsubsection{$C^k$-topology}

The results are proved first using the $C^k$-topology. The Leibniz rule implies that there is a constant $C$ such that for tensors $T$, $T'\in C^k$ 
\begin{equation}\label{Ck continuity of product}
\|T \cdot T' \|_{C^k} \leq C \|T\|_{C^k} \|T'\|_{C^k}.
\end{equation}
The dot stands for any algebraic operation involving tensor product and contraction. The constant $C$ depends only on $k$, not on the metric used to calculated the norms (in contrast to the Sobolev analogue, which is discussed in the next section).

\begin{lemma}\label{Ck continuity of curvature}
There exist positive constants $c$, $K$, such that whenever $g$, $g'$ are two different metrics on the same compact manifold, satisfying
$$
\|g' - g\|_{C^{k+2}} \leq c,
$$
with corresponding curvature tensors $R$, $R'$, then
$$
\| R' - R \|_{C^k} 
\leq 
K \|g' -g\|_{C^{k+2}}.
$$
All norms are taken with respect to the metric $g'$. $K$ depends only on $c$ and $k$ (and not on $g$ or $g'$).
\end{lemma}

\begin{proof}
Let $g =g' + h$. If the corresponding Levi-Civita connections are denoted $\nabla$ and $\nabla'$, then $\nabla = \nabla' +a$, where $a$ corresponds to $-\nabla' h$ under the isomorphism $T^* \otimes \End T \cong T^*\otimes T^* \otimes T^*$ defined by $g$. That is,
$$
a\cdot (g'+h) = -\nabla'h, 
$$
where the dot is an algebraic operation. Hence
\begin{eqnarray*}
\|a\|_{C^{k+1}} 
	& = &
		\|a\cdot g'\|_{C^{k+1}}, \\
	& \leq &
		\| \nabla' h \|_{C^{k+1}} + C \|a\|_{C^{k+1}}\|h\|_{C^{k+1}},
\end{eqnarray*}
(using inequality (\ref{Ck continuity of product}) above). Taking $c <  C^{-1}$ gives
$$
\|a\|_{C^{k+1}} 
\leq 
\|h\|_{C^{k+2}} \left( 1 - c^{-1}\|h\|_{C^{k+2}}\right)^{-1}.
$$

The difference in curvatures is given by $R - R' = \nabla' a + a\wedge a$. The result now follows from (\ref{Ck continuity of product}).
\end{proof}

\begin{lemma}\label{Ck continuity of ricci curvature}
Given $k$ and $M>0$, there exist positive constants $c$ and $K$ such that whenever $g$ and $g'$ are two different metrics on the same compact manifold, satisfying
\begin{eqnarray*}
\|g' - g\|_{C^{k+2}} &\leq& c,\\
\|R'\|_{C^{k}} &\leq& M,
\end{eqnarray*}
where $R'$ is the curvature tensor of $g'$, then
\begin{eqnarray*}
\| \Ric' - \Ric \|_{C^k} 
&\leq &
K \|g' -g\|_{C^{k+2}},\\
\| \Scal' - \Scal \|_{C^k}
&\leq &
K \|g' -g\|_{C^{k+2}}.
\end{eqnarray*}
Here, $\Ric$ and $\Ric'$, $\Scal$ and $\Scal'$ are the Ricci and scalar curvatures of $g$ and $g'$ respectively. 
\end{lemma}

\begin{proof}
The Ricci tensor is given by $\Ric = R \cdot g$ where the dot denotes contraction with the metric. Simple algebra gives
$$
\Ric' -\Ric = (R' -R)\cdot g' - (R' -R) \cdot(g' -g) + R'\cdot(g' -g).
$$
It follows from inequality (\ref{Ck continuity of product}) that $\| \Ric' - \Ric\|_{C^k}$ is controlled by a constant multiple
$$
\|R' - R\|_{C^k} \|g'\|_{C^k} +
\|R' - R\|_{C^k} \|g' -g\|_{C^k} +
\|R'\|_{C^k} \|g' -g\|_{C^k}.
$$
Since the $C^k$-norm of $g'$ is constant, the result follows from Lemma \ref{Ck continuity of curvature}.

A similar argument applies to $\Scal = \Ric \cdot g$.
\end{proof}

\begin{lemma}\label{Ck continuity of derivative}
Given $k$ and $M>0$, there exist positive constants $c$ and $K>0$ such that whenever $(J, \omega)$ and $(J', \omega')$ are two different K\"ahler structures on the same compact manifold satisfying
\begin{eqnarray*}
\|(J', \omega') - (J, \omega)\|_{C^{k+2}} &\leq& c,\\
\|R'\|_{C^{k}}  &\leq& M,
\end{eqnarray*}
where $R'$ is the curvature tensor of $(J', \omega')$, then the linearisations $L$ and $L'$ of the corresponding scalar curvature maps satisfy
$$
\|(L' -L)\phi\|_{L^{p}_{k}} 
\leq 
K\|(J',\omega') - (J, \omega)\|_{C^{k+2}}\|\phi\|_{L^{p}_{k+4}}
$$
for all $\phi \in L^{p}_{k+4}$. All norms are computed with respect to the primed K\"ahler structure.
\end{lemma}

\begin{proof}
The formula (\ref{derivative of scalar curvature}) shows that $L$ is a sum of compositions of the operators $\Delta$, $i\delb\del$, of multiplication by $\Scal$, $\rho$ and $\omega$, and of division of top degree forms by $\omega^n$. It suffices, then, to show that these operations satisfy inequalities analogous to that in the statement of the lemma.  

For multiplication by $\omega$ this is immediate. Since dividing top degree forms by $\omega^n$ is the same as taking the inner product with $\omega^n/n!$ it holds for this operation too. For multiplication by $\Scal$ and $\rho$, the inequalities follow from Lemmas \ref{Ck continuity of ricci curvature} and the inequality $\|u v\|_{L^p_k} \leq C \|u\|_{C^k} \|v\|_{L^p_k}$ for some $C$ (depending only on $k$).

On functions, $\Delta$ is the trace of $i \delb\del$, so to prove the lemma it suffices to prove that $i\delb\del$ satisfies the relevant inequality (\emph{cf.} the proof of Lemma \ref{Ck continuity of ricci curvature}). Since $\pi^{1,0} = \frac{1}{2}(1 - i J)$, the operator $\del = \pi^{1,0} \diff$ satisfies the required inequality, similarly for $\delb$. Hence $i \delb\del$ and $\Delta$ do too. 
\end{proof} 

\subsubsection{$L^p_k$-topology}

In the above discussion of continuity, it is possible to work with Sobolev rather than $C^k$ norms. The same arguments apply with one minor modification. Inequality (\ref{Ck continuity of product}) is replaced by
$$
\|T \cdot T' \|_{L^p_k} \leq C \|T\|_{L^p_k} \|T'\|_{L^p_k},
$$
which holds provided $kp>2n$. Moreover, the constant $C$ depends on the metric through the constants appearing in the Sobolev inequalities
\begin{eqnarray}
\label{C0 sob ineq}
\| S \|_{C^0} &\leq &C' \|S\|_{L^p_k} \quad {\text {for\ }} kp > 2n,\\
\label{Lpk sob ineq}
\| S \|_{L^p_k} &\leq &C''\|S\|_{L^q} \quad {\text {for\ }} kp > 2n.
\end{eqnarray}

With this in mind, the same chain of reasoning which leads to Lemma \ref{Ck continuity of derivative} also proves:

\begin{lemma}\label{Lpk continuity of derivative}
Let $kp - 2n >0$, and $M>0$. There exist positive constants $c$, $K$, such that whenever $(J, \omega)$, $(J', \omega')$ are two different K\"ahler structures on the same compact complex $n$-dimensional manifold, satisfying
\begin{eqnarray*}
\|(J',\omega') - (J, \omega)\|_{L^p_{k+2}}& \leq & c,\\
\|R'\|_{L^{p}_{k}},\, C',\, C'' & \leq & M,
\end{eqnarray*}
where $R'$ is the curvature tensor and $C'$, $C''$ are the Sobolev constants from inequalities (\ref{C0 sob ineq}) and (\ref{Lpk sob ineq}) for the primed K\"ahler structure, then the linearisations $L$ and $L'$ of the corresponding scalar curvature maps satisfy
$$
\|(L' -L)\phi\|_{L^{p}_{k}} 
\leq 
K\|(J',\omega') - (J, \omega)\|_{L^p_{k+2}}\|\phi\|_{L^{p}_{k+4}}
$$
for all $\phi \in L^{p}_{k+4}$. All norms are computed with respect to the primed K\"ahler structure.
\end{lemma}

\section{Approximate solutions}\label{approximate solutions}

Return now to the case where $X$ is a compact connected complex surface and $\pi \colon X \to \Sigma$ is a holomorphic submersion onto a smooth high genus curve with fibres of genus at least $2$.

This section constructs families of metrics on $X$, each depending on a parameter $r$, which have approximately constant scalar curvature. During this procedure, various power series expansions in negative powers of $r$ will be used. Questions of convergence with respect to various Banach space norms will be addressed later. For now, the expression $f(r)=O(r^{-n})$ is to be interpreted as holding pointwise.
When such an expression is used to describe an operator it should be interpreted as holding after the operator acts on a function.

The ultimate aim of this chapter is to construct, for each non-negative integer $n$, a family of metrics $\omega_{r,n}$ parametrised by $r$, satisfying
$$
\Scal(\omega_{r,n}) = -1 +\sum_{i=1}^n c_i r^{-i} + O\left(r^{-n-1}\right),
$$
where $c_1, \ldots, c_n$ are constants. This is accomplished in Theorem \ref{high order approximate solutions}.

\subsection{The first order approximate solution}

Recall the classes
$$
\kappa_r = - 2\pi\left(c_1(V) + r c_1(\Sigma) \right)
$$
mentioned (up to a factor of $2\pi$) in the introduction. Here $V$ denotes the vertical tangent bundle over $X$ and $r$ is a positive real number.

\begin{lemma}\label{kahler classes}
For all sufficiently large $r$, $\kappa_r$ is a K\"ahler class. Moreover, it contains a K\"ahler representative $\omega_r$ whose fibrewise restriction is the canonical hyperbolic metric on that fibre.
\end{lemma}

\begin{proof}
Each fibre has a canonical hyperbolic metric. These metrics define a  Hermitian structure in the holomorphic bundle $V \to X$. Denote the corresponding curvature form by $F_V$, and define a closed real $(1,1)$-form by $\omega_0 = -i F_V$. Notice that $[\omega_0] = -2\pi c_1(V)$, and the fibrewise restriction of $\omega_{0}$ is the hyperbolic metric of that fibre.

Since the fibrewise restriction of $\omega_0$ is non-degenerate, it defines a splitting $TX = V \oplus H$, where
$$
H_x = \{ u \in T_xX : \omega_0(u,v) = 0 \quad \text{for all } v \in V_x \}.
$$
Let $\omega_\Sigma$ be any K\"ahler form on the base, scaled so that $[\omega_\Sigma] = -2\pi c_1(\Sigma)$. The form $\omega_\Sigma$ (pulled back to $X$) is a pointwise basis for the purely horizontal $(1,1)$-forms. This means that, with respect to the vertical-horizontal decomposition, 
$$
\omega_0 = \omega_\sigma \oplus \theta \omega_\Sigma
$$
for some function $\theta \colon X \to \R$, where $\omega_\sigma$ is the hyperbolic K\"ahler form on the fibre $S_\sigma$ over $\sigma$.

For $r > -\inf \theta$, the closed real $(1,1)$-form 
$$
\omega_r = \omega_0 + r\pi^*\omega_\Sigma
$$
is positive, and hence K\"ahler, with $[\omega_r] = \kappa_r$. Its restriction to $S_\sigma$ is $\omega_\sigma$ as required.
\end{proof}

\begin{definition}\label{vertical and horizontal laplacians}
The \emph{vertical Laplacian}, denoted $\Delta_V$, is defined by
$$
\left(\Delta_V \phi \right)\omega_\sigma = i(\delb \del \phi)_{VV},
$$
where $(\alpha)_{VV}$ denotes the purely vertical component of a $(1,1)$-form $\alpha$. The fibrewise restriction of $\Delta_V$ is the Laplacian determined by  $\omega_\sigma$.

The \emph{horizontal Laplacian}, denoted $\Delta_H$, is defined by
$$
\left(\Delta_H \phi \right)\omega_\Sigma = (i\delb \del \phi)_{HH},
$$
where $(\alpha)_{HH}$ denotes the purely horizontal component of a $(1,1)$-form $\alpha$. On functions pulled up from the base $\Delta_H$ is the Laplacian determined by $\omega_\Sigma$.
\end{definition}

\begin{lemma}\label{O(r^-1) expansion of Scal(omega_r)}
\begin{equation}\label{scalar curvature of omega_r}
\Scal (\omega_r)
= 	
-1+ r^{-1}\left(
\Scal(\omega_\Sigma) - \theta + \Delta_V \theta 
\right)
+ O(r^{-2}).
\end{equation}
\end{lemma}

\begin{proof}
The short exact sequence of holomorphic bundles 
$$
0 \to V \to TX \to H \to 0
$$
induces an isomorphism $K_X \cong V^* \otimes H^*$. This implies that the Ricci form of $\omega_r$ is given by $\rho_r = i(F_V +F_H)$, where $F_V$ and $F_H$ are the curvature forms of $V$ and $H$ respectively.

The metric on the horizontal tangent bundle is $(r+\theta)\omega_\Sigma$. Its curvature is given by
$$
iF_H = \rho_\Sigma + i\delb \del \log(1 +r^{-1}\theta),
$$
where $\rho_\Sigma$ is the Ricci form of $\omega_\Sigma$. The curvature of the vertical tangent bundle has already been considered in the definition $\omega_0 = -iF_V$. Hence
\begin{equation}\label{rho_r}
\rho_r 
= 
-\omega_\sigma -\theta\omega_\Sigma + \rho_\Sigma 
+ i\delb \del \log(1 +r^{-1}\theta).
\end{equation}
Taking the trace gives
$$
\Scal(\omega_r) 
= 
-1 + \frac{\Scal(\omega_\Sigma) - \theta}{r + \theta} 
+ \Delta_r \log (1 + r^{-1}\theta).
$$
where $\Delta_r$ is the Laplacian determined by $\omega_r$. Using the formula
\begin{equation}\label{Delta_r}
\Delta_r = \Delta_V + \frac{\Delta_H}{r+\theta}
\end{equation}
and expanding out in powers of $r^{-1}$ proves the result.
\end{proof}

Since $\Scal(\omega_r) = -1 + O(r^{-1})$, setting $\omega_{r,0} = \omega_r$ gives the first family of approximate solutions.

\subsection{The second order approximate solution}

Let $L_r$ denote the linearisation of the scalar curvature map on K\"ahler potentials determined by $\omega_r$. The $r$ dependence of $L_r$ will be of central importance in the proof of Theorem \ref{main theorem}. Its study will essentially occupy the remainder of this article. A first step in this direction is provided by the following lemma. 

\begin{lemma}\label{L_r to O(r^-1)} 
$$
L_r = \Delta_V^2 + \Delta_V + O(r^{-1}).
$$
\end{lemma}

\begin{proof}
Recall the formula (\ref{derivative of scalar curvature}) for $L_r$:
$$
L_r(\phi) = \Delta_r^2\phi -\Scal(\omega_r)\Delta_r\phi + 
\frac{2i\delb \del \phi \wedge \rho_r}{\omega_r^2}.
$$
Equations (\ref{scalar curvature of omega_r}), (\ref{rho_r}) and (\ref{Delta_r}) give the $r$ dependence of $\Scal(\omega_r)$, $\rho_r$ and $\Delta_r$ respectively. Direct calculation gives the result.
\end{proof}

\begin{remark*}
Notice that the $O(1)$ term of $L_r$ is the first order variation in the scalar curvature of the fibres (see Example \ref{L for high genus curves}). This can be seen as an example of the dominance of the local geometry of the fibre in an adiabatic limit. 

Instead of using a calculation as above, this result can be seen directly from formula (\ref{scalar curvature of omega_r}). The $O(1)$ term in $\Scal(\omega_r)$ is $\Scal(\omega_\sigma)$. Rather than considering a K\"ahler potential as a change in $\omega_r$, it can be thought of as a change in $\omega_0$. This gives a corresponding change in $\omega_\sigma$ and the $O(1)$ effect on $\Scal(\omega_r)$ is precisely that claimed. 
\end{remark*}

Given a function $\phi \in C^\infty(X)$, taking the fibrewise mean value gives a function $\pi_\Sigma \phi \in C^\infty(\Sigma)$. The projection maps and $\pi_\Sigma$ and $1-\pi_{\Sigma}$ determine a splitting
$$
C^\infty(X) = C^\infty_0(X) \oplus C^\infty(\Sigma),
$$
where $C^\infty_0(X)$ denotes functions with fibrewise mean value zero.

The previous lemma implies that, at least to $O(r^{-1})$, functions in the image of $L_r$ have fibrewise mean value zero. It is because of this that the $C^\infty_0(X)$ and $C^\infty(\Sigma)$ components of the errors in $\Scal(\omega_r)$ must be dealt with differently.

\subsubsection{The correct choice of $\omega_\Sigma$}
\label{finding omega_Sigma}

Notice that the definition of $\omega_{r}$ has, so far, involved an arbitrary metric on $\Sigma$. 

\begin{theorem}
Each conformal class on $\Sigma$ contains a unique representative $\omega_{\Sigma}$ such that
$$
\pi_{\Sigma} \Scal(\omega_{r}) = -A -r^{-1} + O\left(r^{-2}\right),
$$
where $A$ is the area of a fibre in its hyperbolic metric.
\end{theorem}

\begin{proof}
By Lemma \ref{O(r^-1) expansion of Scal(omega_r)} the $C^{\infty}(\Sigma)$ component of the $O(r^{-1})$ term in $\Scal(\omega_r)$ is $\Scal(\omega_\Sigma) - \pi_\Sigma \theta$.
It follows from the definition of $\theta$ (as the horizontal part of $\omega_0 = -iF_V$ divided by $\omega_\Sigma$) that 
$$
\pi_\Sigma \theta = - A^{-1}\Lambda_\Sigma \pi_*(F_V^2),
$$
where $\Lambda_\Sigma$ is the trace on $(1,1)$-forms on $\Sigma$ determined by $\omega_\Sigma$. 

Write $\alpha = -\pi_{*}(F_{V}^{2})$
The surface $X$ determines a map to the moduli space of curves. The form $\alpha$ is a representative for the pull back of the first tautological class via this map. (See, \emph{e.g.}\ \cite{harris.morrison:moc}.) Since the first tautological class is ample, $\int \alpha\geq 0$. The results of \cite{kazdan.warner:icfpdewatrg} (which discusses prescribing curvature on curves) can now be applied. They prove the existence of a unique metric in each conformal class with $\Scal - A^{-1}\Lambda \alpha = -1$.
\end{proof}

From now on this choice of metric is assumed to be included in the definition of $\omega_r$.

\subsubsection{The correct choice of K\"ahler potential $\phi_1$}

Let $\Theta_1$ denote the $C^\infty_0(X)$ component of the $O(r^{-1})$ term in $\Scal(\omega_r)$. This means that
$$
\Scal(\omega_r) = -1 +r^{-1}(\Theta_1 -1) +O(r^{-2}).
$$
It follows from Lemma \ref{L_r to O(r^-1)} that
\begin{equation}\label{Scal(omega_r) to O(r^-2)}
\Scal(\omega_r + ir^{-1}\delb \del \phi) 
= 
\Scal(\omega_r) +r^{-1}(\Delta_V^2 + \Delta_V)\phi + O(r^{-2}).
\end{equation}

\begin{lemma}\label{invertibility of O(r^-1) piece of L_r}
Let $\Theta \in C^\infty_0(X)$. There exists a unique $\phi \in C^\infty_0(X)$ such that
$$
\left(\Delta_V^2 + \Delta_V\right) \phi = \Theta.
$$
\end{lemma}

\begin{proof}
Given a function $\phi \in C^\infty(X)$, let $\phi_\sigma$ denote the restriction of $\phi$ to $S_\sigma$. The fibrewise restriction of the operator $\Delta_V^2+\Delta_V$ is the first order variation of the scalar curvature of the fibre. Applying Example \ref{L for high genus curves} fibrewise certainly gives a unique function $\phi$ on $X$ such that $\phi$ has fibrewise mean value zero; for each $\sigma$, $\phi_\sigma \in C^\infty(S_\sigma)$ and $L_\sigma \phi_\sigma = \Theta_\sigma$, \emph{i.e.}\ $(\Delta_V^2 + \Delta_V)\phi = \Theta$. It only remains to check that $\phi$ is smooth transverse to the fibres. (The operator $\Delta_V^2 + \Delta_V$ is only elliptic in the fibre directions, so regularity only follows automatically in those directions.)

In fact, this is straight forward. Since $\phi_\sigma = L^{-1}_\sigma \Theta_\sigma$ the required differentiability follows from that of $\Theta$ and the fact that $L_\sigma$ is a smooth family of differential operators.
\end{proof}

Applying this lemma to $\Theta = -\Theta_1$ and using equation (\ref{Scal(omega_r) to O(r^-2)}) shows that there exists a unique $\phi_1 \in C^\infty_0(X)$ such that
the metric $\omega_{r,1} = \omega_r + i\delb\del r^{-1}\phi_1$ is an $O(r^{-2})$ approximate solution to the constant scalar curvature equation:
$$
\Scal(\omega_{r,1}) = -1 -r^{-1} + O(r^{-2}).
$$

\subsection{The third order approximate solution}\label{third}

Now that the correct metric has been found on the base, the higher order approximate solutions are constructed recursively. In order to demonstrate the key points clearly, however, this section does the first step in detail. 

The strategy is straightforward, even if the notation sometimes isn't. The first step is to find a K\"ahler potential $f_1$ on the base to deal with the $C^\infty(\Sigma)$ component of the $O(r^{-2})$ error; that is, so that
$$
\Scal(\omega_{r,1} + i\delb\del f_1)
=
-1 -r^{-1} + (c +\Theta'_2)r^{-2} + O(r^{-3}),
$$
for some constant $c$, where $\Theta'_2$ has fibrewise mean value zero. 

The second step is to find a K\"ahler potential $\phi_2$ to deal with the remaining $O(r^{-2})$ error $\Theta'_2$; that is, so that
$$
\Scal\left(\omega_{r,1} + i\delb\del (f_1 + r^{-2}\phi_2)\right)
=
-1 - r^{-1} + cr^{-2} + O(r^{-3}).
$$

Both of the potentials $f_1$ and $\phi_2$ are found as solutions to linear partial differential equations. To find the relevant equations, it is important to understand the linearisation of the scalar curvature map on K\"ahler potentials determined by $\omega_{r,1}$ (and the operators determined by the later, higher order, approximate solutions). To this end, the first lemma in this section deals with the $r$ dependence of such an operator when the fibrewise metrics are not necessarily the canonical constant curvature ones. First, some notation.

\subsubsection*{Notation for Lemma \ref{expansion of L(Omega_r)}}

Let $\Omega_0$ be any closed real $(1,1)$-form whose fibrewise restriction is K\"ahler. Let $\Omega_\sigma$ be the K\"ahler form on $S_\sigma$ induced by $\Omega_0$. Let $\Omega_\Sigma$ be any choice of metric on the base. As in Lemma \ref{kahler classes}, for large enough $r$, the form $\Omega_r = \Omega_0 + r\Omega_\Sigma$ is K\"ahler; the vertical-horizontal decomposition of the tangent bundle determined by $\Omega_r$ depends only on $\Omega_0$.

\begin{definition}
The form $\Omega_\Sigma$ is a pointwise basis for the horizontal $(1,1)$-forms. Define a function $\xi$ as follows: write the vertical-horizontal decomposition of $\Omega_0$ (with respect to $\Omega_r$) as $\Omega_0 = \Omega_\sigma \oplus \xi \Omega_\Sigma$.
\end{definition}

\begin{definition}
The family of fibrewise K\"ahler metrics $\Omega_\sigma$ determines a Hermitian structure in the vertical tangent bundle. Denote the curvature of this bundle as $F_V$. Define a function $\eta$ as follows: write the horizontal-horizontal component of $iF_V$ (with respect to $\Omega_r$) as $\eta \Omega_\Sigma$.
\end{definition}

\begin{remark*}
Since the fibrewise metrics are not the canonical constant curvature ones, this curvature form is not the same as that appearing earlier. If, instead of any old $\Omega_0$ and $\Omega_\Sigma$, the forms $\omega_0$ and $\omega_\Sigma$ from before are used in both of these definitions, then $\xi = -\eta = \theta $.
\end{remark*}

\begin{definition}
Taking the fibrewise mean value of $\eta$ gives a function $\pi_\Sigma \eta$ on the base. Using this, define a fourth order differential operator
$$
D_\Sigma : C^\infty(\Sigma) \to C^\infty(\Sigma), 
$$
$$
D_\Sigma(f) 
= 
\Delta^2_\Sigma f - (\Scal(\Omega_\Sigma) +\pi_\Sigma \eta)\Delta_\Sigma f,
$$
where $\Delta_\Sigma$ is the $\Omega_\Sigma$-Laplacian.
\end{definition}

\begin{remark}\label{D_Sigma is invertible}
The operator $D_\Sigma$ is the linearisation of a nonlinear map on functions, which is now described. Taking the fibrewise mean value of $\eta\Omega_{\Sigma}$ defines a $2$-form on the base $\Sigma$. Notice this is independent of the choice of $\Omega_\Sigma$. The trace of this form with respect to $\Omega_\Sigma$ is precisely the fibrewise mean value of $\eta$. 

Next, consider varying $\Omega_\Sigma$ by a K\"ahler potential $f \in C^\infty(\Sigma)$. Denote by $\Lambda_{\Sigma, f}$ the trace operator determined by $\Omega_\Sigma + i\delb\del f$. The equation
$$
\Lambda_{\Sigma,f} = \frac{\Lambda_{\Sigma}}{1 + \Delta_\Sigma f}
$$
shows that the linearisation at $0$ of the map $f \mapsto \pi_\Sigma \eta$ is $-\pi_\Sigma \eta \Delta_\Sigma$. Combining this with the formula for the linearisation of the scalar curvature map on curves derived in Example \ref{L for high genus curves}, shows that $D_\Sigma$ is the linearisation, at $0$, of the map
$$
F \colon f \mapsto \Scal(\Omega_\Sigma + i\delb\del f) + \pi_\Sigma \eta.
$$

If, instead of using any old $\Omega_0$, the definition were made using $\omega_0$ from earlier then the map $F$ is one which has been described before. It is precisely the map which was shown to take the value $-1$ at $\omega_\Sigma$ (see section \ref{finding omega_Sigma}). Notice that using $\omega_0$ and $\omega_\Sigma$ to define $D_\Sigma$ gives $D_\Sigma = \Delta_\Sigma^2 + \Delta_\Sigma$. As in Example \ref{L for high genus curves}, this operator is an isomorphism on functions of mean value zero (when considered as a map between the relevant Sobolev spaces).
\end{remark}

The vertical and horizontal Laplacians are defined just as before, with $\Omega_0$ and $\Omega_\Sigma$ replacing $\omega_0$ and $\omega_\Sigma$ respectively (see Definition \ref{vertical and horizontal laplacians}). To indicate that they are defined with respect to different forms (and also a different vertical-horizontal decomposition of the tangent bundle, notice), the vertical and horizontal Laplacians determined by $\Omega_0$ and $\Omega_\Sigma$ are denoted $\Delta_V'$ and $\Delta_H'$. The un-primed symbols are reserved for the vertical and horizontal Laplacians determined by $\omega_0$ and $\omega_\Sigma$. Let $L(\Omega_r)$ denote the linearisation of the scalar curvature map on K\"ahler potentials defined by $\Omega_r$. 

\begin{lemma}\label{expansion of L(Omega_r)}
$$
L(\Omega_r) = (\Delta_V^{'2} - \Scal(\Omega_\sigma) \Delta'_V) + r^{-1}D_1 + r^{-2}D_2 +O(r^{-3}),
$$
where the operators $D_1$ and $D_2$ have the following behaviour: if $f$ is a function pulled back from $\Sigma$,
\begin{eqnarray}
\label{D_1 for Omega_r}
D_1(f)	&	=	&	0,\\
\label{D_2 for Omega_r}
\pi_\Sigma D_2(f)	&	=	&	D_\Sigma(f).
\end{eqnarray}
\end{lemma}

\begin{proof}
The proof given here is a long calculation. A slightly more conceptual proof is described in a following remark. Recall the formula (\ref{derivative of scalar curvature}) for the linearisation of the scalar curvature map. It involves the Laplacian, the scalar curvature and the Ricci form of $\Omega_r$. Repeating the calculations that were used when the fibres had cscK metrics gives formulae for these objects. They are, respectively,
\begin{eqnarray}
\label{laplacian of Omega_r}
\Delta_{\Omega_r} & = & \Delta_V' 
+ 
\frac{\Delta_H'}{r + \xi}, \\
\label{scalar curvature of Omega_r}
\Scal(\Omega_r) & = & 
\Scal(\Omega_\sigma) 
+ 
\frac{\Scal(\Omega_\Sigma) + \eta}{r +\xi}
+
\Delta_{\Omega_r} \log (1 + r^{-1} \xi), \\
\label{ricci form of Omega_r}
\rho(\Omega_r) & = & \rho(\Omega_\sigma) 
+
\left(\Scal(\Omega_\Sigma) +\eta\right)\Omega_\Sigma
+ i\delb\del \log(1 + r^{-1}\xi).
\end{eqnarray}

The result now follows from routine manipulation and expansion of power series; the following formulae can be verified for $D_1$ and $D_2$:
\begin{eqnarray*}
D_1 & = & 2\Delta'_V\Delta'_H - (\Delta'_V \xi) \Delta'_V,\\
D_2 & = & \Delta_H^{'2} -(\Scal(\Omega_\Sigma) + \eta)\Delta'_H - \eta\Delta'_V\Delta'_H \\
	  &    & + \frac{1}{2}(\Delta'_V (\xi^2)) \Delta'_V + (\Delta'_V \eta) \Delta'_H.
\end{eqnarray*}
The statements about $D_1(f)$ and $\pi_\Sigma D_2(f)$ for $f$ pulled up from the base follow from these equations.
\end{proof}

\begin{remark*}
The actual equations for $D_1$ and $D_2$ will not be needed in what follows. All that will be used is their stated behaviour on functions on the base as stated in  Lemma \ref{expansion of L(Omega_r)}. This behaviour can be understood, without laborious calculation, as follows. 

The potential $f$ can be thought of as altering $\Omega_{\Sigma}$, rather than $\Omega_{r}$. Since $\Omega_{\Sigma}$ is scaled by $r$ in the definition of $\Omega_{r}$, the effect is equivalent to adding the potential $r^{-1}f$ to $\Omega_{\Sigma}$. The analogue of equation (\ref{scalar curvature of omega_r}) shows that the lowest order effect of $\Omega_{\Sigma}$ on $\Scal(\Omega_{r})$ occurs at $O(r^{-1})$. Hence the combined effect is $O(r^{-2})$: for potentials $f$ pulled up from the base $D_1(f)=0$. 

The fibrewise mean value of the $O(r^{-1})$ term in $\Scal(\Omega_r)$ is 
$$\Scal(\Omega_\Sigma) +\pi_\Sigma \eta.$$ 
So, after taking the fibrewise mean value, a change of $r^{-1}f$ in $\Omega_\Sigma$ gives a change in $\pi_\Sigma \Scal(\Omega_r)$ whose $O(r^{-2})$ term is given by the derivative of the above expression with respect to K\"ahler potentials on the base. Hence $\pi_\Sigma D_2 (f)= D_\Sigma(f)$.
\end{remark*}

\subsubsection{The correct choice of K\"ahler potential $f_1$}

Let $L_{r,1}$ be the linearisation of the scalar curvature map on K\"ahler potentials determined by $\omega_{r,1}$.

\begin{lemma}\label{O(r^-2) expansion of L_r,1}
Let $f \in C^\infty(\Sigma)$. Then
$$
\pi_\Sigma L_{r,1} (f) = r^{-2}(\Delta^2_\Sigma + \Delta_\Sigma)f + O(r^{-3}).
$$
\end{lemma}

\begin{proof}
Begin by applying Lemma \ref{expansion of L(Omega_r)} with 
\begin{eqnarray*}
\Omega_0			&	=	&	\omega_0 + i\delb\del r^{-1}\phi_1,\\
\Omega_\Sigma	&	=	&	\omega_\Sigma.
\end{eqnarray*}
There is a slight difficultly in interpreting the expansion given in Lemma \ref{expansion of L(Omega_r)}. The $r$-dependence of $\Omega_0$ means that the coefficients in the $O(r^{-3})$ piece of that expansion will be $r$-dependent, \emph{a priori} making them of higher order overall. 

In fact, this can't happen. The reason is that all such coefficients come ultimately from analytic expressions in the fibrewise metrics induced by $\Omega_0$ (as is shown, for example, by the calculation described in the proof of Lemma \ref{expansion of L(Omega_r)}). These metrics have the form
$$
\Omega_\sigma = (1 + r^{-1}\Delta_V \phi_1)\omega_\sigma.
$$
(Here $\Delta_V$ is the vertical Laplacian determined by $\omega_0$.) Since the fibrewise metric is algebraic in $r^{-1}$, the coefficients in the expression form Lemma \ref{expansion of L(Omega_r)} are analytic in $r^{-1}$, \emph{i.e.}\ they have expansions involving only non-positive powers of $r$.

This means that the $O(r^{-2})$ term can simply be read off from the formula given in Lemma \ref{expansion of L(Omega_r)}. This gives
$$
\pi_\Sigma L_{r,1}(f) = r^{-2} D_\Sigma(f) + O(r^{-3}).
$$
As is pointed out in Remark \ref{D_Sigma is invertible}, for the choice of $\omega_\Sigma$ that was determined whilst finding the $O(r^{-2})$ approximate solution, $D_\Sigma = \Delta^2_\Sigma + \Delta_\Sigma$ as required.
\end{proof}

Denote the $C^\infty(\Sigma)$ component of the $O(r^{-2})$ term of $\Scal(\omega_{r,1})$ by $\Theta_2$:
$$
\pi_\Sigma \Scal(\omega_{r,1}) = -1 - r^{-1} + r^{-2}\Theta_2 +O(r^{-3}).
$$
Let $c$ be the mean value of $\Theta_{2}$ and let $f_{1}$ solve $(\Delta_{\Sigma}^{2}+\Delta_{\Sigma})f_{1}= c - \Theta_{2}$. (Existence of $f_{1}$ follows from  Example \ref{L for high genus curves}.) Elliptic regularity ensures that $f_1$ is smooth, completing the first step in finding the $O(r^{-3})$ approximate solution: by Lemma \ref{O(r^-2) expansion of L_r,1},
$$
\Scal(\omega_{r,1} + i\delb\del f_1) 
= 
-1-r^{-1}+r^{-2}(c + \Theta'_2) +O(r^{-3}),
$$
where $\Theta'_2$ has fibrewise mean value zero.

\subsubsection{The correct choice of K\"ahler potential $\phi_2$}

Let $L'_{r,1}$ be the linearisation of the scalar curvature map on K\"ahler potentials determined by the metric $\omega_{r,1} + i\delb\del f_1$.

\begin{lemma}\label{O(r^-1) expansion of L'_r,1}
$$
L'_{r,1} = \Delta_V^2 + \Delta_V + O(r^{-1}).
$$
\end{lemma}

\begin{remark*}
Again, the symbol $\Delta_V$ means the vertical Laplacian determined by the form $\omega_0$. Compare this with Lemma \ref{L_r to O(r^-1)}.
\end{remark*}

\begin{proof}
Apply Lemma \ref{expansion of L(Omega_r)} with 
\begin{eqnarray*}
\Omega_0			&	=	&	\omega_0 + i\delb\del r^{-1}\phi_1,\\
\Omega_\Sigma	&	=	&	\omega_\Sigma + i\delb\del r^{-1}f_1.
\end{eqnarray*}
As in the proof of Lemma \ref{O(r^-2) expansion of L_r,1}, there is a problem with interpreting the expansion in Lemma \ref{expansion of L(Omega_r)}, namely that the $r$-dependence of $\Omega_0$ and $\Omega_\Sigma$ means that the coefficients in the expansion are also $r$-dependent. As in the proof of Lemma \ref{O(r^-2) expansion of L_r,1}, however, this actually causes no difficulty. Both forms are algebraic in $r^{-1}$, hence the coefficients in the expansion are analytic in $r^{-1}$; the $r$-dependence of the coefficient of $r^{-n}$ causes changes only at $O(r^{-n-k})$ for $k\geq 0$. This means that the genuine $O(1)$ behaviour of $L'_{r,1}$ is the same as the $O(1)$ behaviour of $\Delta_V^{'2} - \Scal(\Omega_\sigma) \Delta'_V$. Here, $\Omega_\sigma$ is the metric on $S_\sigma$ determined by $\Omega_0$, \emph{i.e.}\ $\Omega_\sigma = (1 + r^{-1}\Delta_V \phi_1)\omega_\sigma$. 

Since, to $O(1)$, $\Omega_\sigma$ and $\omega_\sigma$ agree, 
\begin{eqnarray*}
\Scal(\Omega_\sigma)	&	=	& \Scal(\omega_\sigma) +O(r^{-1}),\\
\Delta_V'				&	=	& \Delta_V + O(r^{-1}).
\end{eqnarray*}
Hence
$$
\Delta_V^{'2} - \Scal(\Omega_\sigma) \Delta'_V
=
\Delta^2_V + \Delta_V + O(r^{-1}).
$$
\end{proof}

Let $\omega_{r,2} = \omega_{r,1} + i\delb\del(f_1+ r^{-2}\phi_2)$, where $\phi_{2}$ solves $(\Delta_{V}^{2}+\Delta_{V})\phi_{2} = -\Theta'$. (Existence of $\phi_{2}$ follows from Lemma \ref{invertibility of O(r^-1) piece of L_r}.) Then
$$
\Scal(\omega_{r,2}) 
= 
-1 -r^{-1} + c r^{-2} + O(r^{-3}).
$$

\subsection{The higher order approximate solutions}\label{higher}

\begin{theorem}[Approximately cscK metrics]
\label{high order approximate solutions}
Let $n$ be a positive integer. There exist functions $f_1, \ldots , f_{n-1} \in C^\infty(\Sigma)$, $\phi_1, \ldots ,\phi_n \in C^\infty_{0}(X)$ and constants $c_{1}, \dots c_{n}$ such that the metric
$$
\omega_{r,n} = \omega_r + i\delb\del\sum_{i=1}^{n-1}r^{-i+1}f_i
+ i\delb\del\sum_{i=1}^n r^{-i}\phi_i
$$
satisfies
$$
\Scal(\omega_{r,n}) = -1 + \sum_{i=1}^{n}c_i r^{-i} + 
O\left(r^{-n-1}\right).
$$
\end{theorem}

\begin{proof}
The proof is by induction, the inductive step being the same as that used to construct the third order approximate solution above. 
\end{proof}

In fact, it is straight forward to show that the functions $f_i$, $\phi_i$ are unique subject to the constraints
$$
\int_\Sigma f_i\,\omega_\Sigma = 0
\qquad
\int_{S_\sigma} \phi_i\,\omega_\sigma = 0.
$$
This uses the injectivity of the operator discussed in Example \ref{L for high genus curves} (acting on functions with mean value zero).

It is also possible to calculate the exact values of the $c_{i}$ by considering the mean value of $\Scal(\omega_{r})$. Let 
$
\vol_{r} =\int \omega^{2}_{r}/2 = rA + B
$
where $A=[\omega_{0}].[\omega_{\Sigma}]$ and $B = \frac{1}{2}[\omega_{0}]^{2}$, and let
$
\int \Scal(\omega_{r})\, \omega_{r}^{2}/2 
= 
2\pi c_{1}(X).[\omega_{r}] 
= 
rC + D
$
where $C =2\pi c_{1}(X).[\omega_{0}]$ and $D = 2\pi c_{1}(X).[\omega_{\Sigma}]$. Then the mean value of the scalar curvature is $(rC+D)(rA+B)^{-1}=-1 + \sum c_{i}r^{-i}$ where 
$$
c_{i}=(-1)^{i}A^{-i}B^{i}(A^{-1}C - B^{-1}D).
$$

\subsection{Summary}\label{summary of approximate solutions}

Four essential facts were used in the construction of the approximate solutions in Theorem \ref{high order approximate solutions}.
\begin{enumerate}
\item The nonlinear partial differential equation $\Scal(\omega_\sigma) = \text{const.}$ in the fibre directions has a solution. This enables $\omega_{r,0}$ to be constructed. 
\item The linearisation of this equation, at a solution, is surjective onto functions with mean value zero. This enables the $C^\infty_0(X)$ components of error terms to be eliminated. 
\item The nonlinear partial differential equation $\Scal(\omega_\Sigma) - \Lambda_\Sigma \alpha =\text{const.}$
on the base has a solution. This enables $\omega_{r,1}$ to be constructed. 
\item The linearisation of this equation, at a solution, is surjective onto functions with mean value zero. This enables the $C^\infty(\Sigma)$ components of error terms to be made constant.
\end{enumerate}

The equation in (3) and (4) is needed only because the linearisation of (2) does not map onto functions pulled up from the base. The surjectivity of the linear operators in (2) and (4) can be viewed in terms of automorphisms of the solutions in (1) and (3) respectively. This is because both operators are elliptic with index zero; they are surjective if and only if they have no kernel (thought of as maps between spaces of functions with mean value zero). The absence of any kernel is equivalent to there being no nontrivial family of solutions to the equations in (1) and (3).

Finally, it should be noted that the two parameters $r$ and $n$ appearing in the approximate solutions are of a very different nature. In particular, whilst the perturbation to a genuine solution is carried out, $n$ is considered as fixed, whilst $r$ tends to infinity.

\section{Applying the inverse function theorem}
\label{set-up ift}

First, some notation. The integer $n$ is considered as fixed throughout and so is often omitted. Write $g_r$ for the metric tensor corresponding to the K\"ahler form $\omega_{r,n}$. Each metric defines Sobolev spaces $L^2_k(g_r)$ of functions over $X$. Since the Sobolev norms determined by $g_r$, for different values of $r$, are equivalent, the spaces $L^2_k(g_r)$ contain the same functions. The constants of equivalence, however, will depend on $r$. Similarly for the Banach spaces $C^k(g_r)$. When the actual norms themselves are not important, explicit reference to the metrics will be dropped, and the spaces referred to simply as $C^k$ or $L^2_k$.

Statements like ``$a_r \to 0$ in $C^k(g_r)$ as $r \to \infty$,'' mean ``$\|a_r\|_{C^k(g_r)} \to 0$ as $r \to \infty$.'' Notice that both the norm and the object whose norm is being measured are changing with $r$. Similarly statements such as ``$a_r$ is $O(r^{-1})$ in $L^2_k(g_r)$ as $r \to \infty$,'' mean ``$\|a_r\|_{L^2_k(g_r)}$ is $O(r^{-1})$ as $r \to \infty$.''

Let $L_{r}$ denote the linearisation of the scalar curvature map on K\"ahler potentials determined by $\omega_{r,n}$. As with Example \ref{L for high genus curves}, this derivative will be shown to be an isomorphism when considered modulo the constant functions.

Recall from Remark \ref{remarks on L} that when the scalar curvature of a K\"ahler metric $\omega$ is constant, the corresponding linear operator maps into functions with $\omega$-mean value zero. In the situation considered here, $\omega_{r,n}$ has nearly constant scalar curvature. So it makes sense to try and show that $L_r$ gives an isomorphism after composing with the projection $p$ onto functions with $\omega_{r,n}$-mean value zero.

Let $L^2_{k,0}$ denote functions in $L^2_k$ with $\omega_{r,n}$-mean value zero. Composing the scalar curvature map with the projection $p$ gives, for $k>1$, a map
$$
S_r \colon L^2_{k+4,0}\to L^2_{k,0},
\quad
S_r(\phi) = p\Scal(\omega_{r,n} + i\delb\del \phi).
$$
To complete the proof of Theorem \ref{main theorem} it will be shown that for each $k$ and sufficiently large $r$ there is a unique $\phi\in L^2_{k+4,0}$ with $S_r(\phi)=0$.

In order to apply the inverse function theorem to find $\phi$ it is necessary to know that $S_{r}(0)$ is sufficiently close to zero. Assume for now that $S_{r}(0) \to 0$ in $L^{2}_{k}(g_{r})$. (Section \ref{approximate solutions} shows only that $S_{r}(0)$ converges to zero pointwise. Convergence in $L^{2}_{k}(g_{r})$ is proved in Lemma \ref{approximate solution in Banach space}.) Assume also that $pL_{r}$, the linearisation of $S_{r}$, is an isomorphism between spaces of functions with mean value zero. (This is proved in Theorem \ref{bounded inverse}.) These two facts alone, however, are not sufficient to be able to apply the inverse function theorem. A close look at the statement of the inverse function theorem may clarify why.

\begin{theorem}[Quantitative inverse function theorem]
\label{ift}
~
\begin{itemize}
\item
Let $F:B_1 \to B_2$ be a differentiable map of Banach spaces, whose derivative at $0$, $DF$, is an isomorphism of Banach spaces, with inverse $P$. 
\item
Let $\delta'$ be the radius of the closed ball in $B_1$, centred at $0$, on which $F-DF$ is Lipschitz, with constant $1/(2\|P\|)$.
\item
Let $\delta =\delta'/(2\|P\|)$. 
\end{itemize}
Then whenever $y\in B_2$ satisfies $\|y-F(0)\| < \delta$, there exists $x \in B_1$ with $F(x)=y$. Moreover such an $x$ is unique subject to the constraint $\|x\| < \delta'$.
\end{theorem}

(This quantitative statement of the inverse function theorem follows from the standard proof; see, for example, \cite{schechter:hoaaif}.)

Applying this to the map $S_{r}\colon L^{2}_{k+4,0} \to L^{2}_{k,0}$ (and assuming its derivative is an isomorphism) gives the existence of a $\delta_{r}$ such that if $\|S_{r}(0)\|_{L^2_k(g_{r})} < \delta_{r}$, 
then there exists a $\phi$ with $S_{r}(\phi)=0$. The proof will be completed by showing that (for any choice of $n \geq 6$) $S_{r}(0)$ converges to zero more quickly than $\delta_{r}$. 

\section{Local analysis}\label{local analysis}

To control the constants appearing in the local analytic estimates, this section constructs a local (over the base) model for the metrics $\omega_{r,n}$.

The notation used here is all defined in Section \ref{approximate solutions}. In particular the forms $\omega_0$, $\omega_r$ and the function $\theta$ are defined in the proof of Lemma \ref{kahler classes}, whilst the higher order approximate solutions $\omega_{r,n}$ are constructed in Theorem \ref{high order approximate solutions}.

\subsection{Constructing the local model}

Let $D \subset \Sigma$ be a holomorphic disc centred at $\sigma_0$. Since $D$ is contractible, $X|_D$ is diffeomorphic to $S \times D$. The horizontal distribution of $\omega_0$ is trivial when restricted to the central fibre $S_{\sigma_0}$. By applying a further diffeomorphism if necessary the identification $X|_D \cong S \times D$ can be arranged so that the horizontal distribution on $S_{\sigma_0}$ coincides with the restriction to $S_{\sigma_0}$ of the $TD$ summand in the splitting
\begin{equation}\label{flat splitting}
T(S\times D) \cong TS \oplus TD.
\end{equation}

For each value of $r$, there are two K\"ahler structures on $S \times D$ of interest. The first comes from simply restricting the K\"ahler structure $(X, J, \omega_{r,n})$ to $X|_D$. The complex structure has the form $J = J_\sigma \oplus J_D$ with respect to the splitting (\ref{flat splitting}), where $J_D$ is the complex structure on $D$ and $J_\sigma$ is the varying complex structure on the fibres.

The second is the natural product structure. With respect to (\ref{flat splitting}), let
\begin{eqnarray*}
J' & = & J_S \oplus J_D, \\
\omega_r'& = & \omega_S \oplus r\omega_D.
\end{eqnarray*}
where $\omega_D$ is the flat K\"ahler form on $D$ agreeing with $\omega_\Sigma$ at the origin, and $J_S$, $\omega_S$ are the complex structure and K\"ahler form on the central fibre $S = S_{\sigma_0}$. Denote by $g'_r$ the corresponding metric on $S\times D$.

It is useful to have the following result stated explicitly.

\begin{lemma}\label{behaviour of C^k(g'_r)}
Let $\alpha \in C^{k}(T^{*}X^{\otimes i})$. Over $X|_{D}$, $\|\alpha\|_{C^k(g'_r)} = O(1)$. Moreover, if $\alpha$ is pulled up from the base, $\|\alpha\|_{C^k(g'_r)} = O\left(r^{-i/2}\right)$.
\end{lemma}

\begin{proof}
These statements follow from the fact that $g'_r$ is a product metric, scaled by $r$ in the $D$ directions.
\end{proof}

\begin{theorem}\label{local model theorem}
For all $\varepsilon >0$, $\sigma_0 \in \Sigma$, there exists a holomorphic disc $D \subset \Sigma$, centred at $\sigma_0$, such that for all sufficiently large $r$, over $X|_D$,
$$
\|(J', \omega'_r) - (J, \omega_{r,n})\|_{C^k(g'_r)} < \varepsilon.
$$
\end{theorem}

\begin{proof}
First notice that, by Lemma \ref{behaviour of C^k(g'_r)}, for any holomorphic disc $D \subset \Sigma$, over $X|_D$, 
$$
\|\omega_{r,n} - \omega_{r,0}\|_{C^k(g'_r)} = O(r^{-1}).
$$
Since $\omega_{r,n} - \omega_{r,0}$ is $O(r^{-1})$ in $C^k(g'_r)$, it suffices to prove the result just for $n=0$. 

Choose a holomorphic disc $D$ centred at $\sigma_0$. The splitting (\ref{flat splitting}) is parallel, implying that 
$$
\nabla^i(J-J') \in \End(TS) \otimes T^*(S \times D)^{\otimes i}.
$$
The only changes in length as $r$ varies come from the $T^*$ factor. Write $\nabla^i(J - J')$ as  $\alpha_i + \beta_i$ with respect to the splitting
$$ 
T^*(S\times D)^{\otimes i} 
\cong 
T^*S^{\otimes i} 
	\oplus 
\Big(
	\left(T^*S^{\otimes i-1}\otimes T^*D \right) 
		\oplus \cdots \oplus 
	T^*D^{\otimes i} 
\Big).
$$
$$
\alpha_i \in T^*S^{\otimes i}, \quad \quad 
\beta_i \in 
	\left(T^*S^{\otimes i-1}\otimes T^*D \right) 
		\oplus \cdots \oplus 
	T^*D^{\otimes i}.
$$

The metric $g'_r$ does not change in the $S$-directions, so $|\alpha_i|_{g'_r}$ is independent of $r$. Since $J=J'$ on the central fibre, reducing the size of $D$ ensures that $|\alpha_i|_{g'_r}$ is less than $\varepsilon /(2k+2)$.

The metric $g'_r$ scales lengths of cotangent vectors by $r^{-1/2}$ in the base directions. So $|\beta_i|_{g'_r} =O(r^{-1/2})$; for large enough $r$, $|\beta_i|_{g'_r}$ is less than $\varepsilon/(2k+2)$. Hence 
$$
\|\nabla^i(J-J')\|_{C^0(g'_r)} < \varepsilon/(k+1).
$$ 
Summing from $i=0, \ldots k$ proves $\|J' - J\|_{C^k(g'_r)} < \varepsilon$.

To prove $\|\omega'_r -\omega_{r}\|_{C^k(g'_r)} < \varepsilon$ it is enough to prove the same result for the metrics $g_{r}$, $g'_r$ (since the K\"ahler forms can be recovered algebraically from the metric tensors via the complex structures).

Let $u_1, u_2$ be a local $g_S$-orthonormal frame for $TS$ and $v_1, v_2$ be a local $g_D$-orthonormal frame for $TD$. Recall $g_r$ induces a different horizontal-vertical splitting of the tangent bundle of $S \times D$, which is independent of $r$. With respect to this splitting, $g_{r} = g_\sigma \oplus (r + \theta)g_\Sigma$, where $g_\sigma$ is the hyperbolic metric on the fibre $S_\sigma$, and $g_\Sigma$ the metric on the base. Write $v_j = \eta_j + \xi_j$ with respect to the horizontal-vertical splitting induced by $g_r$ ($\eta_j$ is horizontal, $\xi_j$ is vertical). With respect to the $g'_r$-orthonormal frame $u_1, u_2, r^{-1/2}v_1, r^{-1/2}v_2$, the matrix representative for $g_{r}$ is
$$
\left(
\begin{array}{ccc}
g_\sigma (u_i, u_j) && r^{-1/2}g_\sigma (u_i, \xi_j) \\ && \\
r^{-1/2}g_\sigma(u_j, \xi_i) && 
\left( 1 + r^{-1}\theta \right) g_\Sigma(\eta_i, \eta_j )
+ r ^{-1} g_\sigma(\xi_i, \xi_j )
\end{array}
\right)
$$

This means that, in a $g'_r$-orthonormal frame, $g_{r} - g'_r$ has the matrix representative
$$
\left(
\begin{array}{cc}
g_\sigma(u_i, u_j) - \delta_{ij} & 0 \\
 & \\
0 & g_\Sigma(\eta_i, \eta_j ) - \delta_{ij}
\end{array}
\right)
+ r^{-1/2}A + r^{-1}B.
$$
for fixed matrices $A$ and $B$. 

The top left corner of the first term vanishes along the central fibre. Just as in the proof of $\|J' -J\|<\varepsilon$,  this can be made arbitrarily small in $C^k(g'_r)$ by shrinking $D$ and taking $r$ large. 

The bottom right corner of the first term is a function of the $D$-variables only. By construction it vanishes at the origin. The $C^0(g'_r)$-norm of this piece is just the conventional $C^0$-norm of the function $g_\Sigma(\eta_i, \eta_j) - \delta_{ij}$ and hence can be made arbitrarily small by shrinking $D$.  

The derivatives of this piece are all in the $D$-directions. The length of the $i$-th derivative is $O(r^{-i/2})$ due to the scaling of $g'_r$ in the $D$-directions. Hence the $C^k(g'_r)$ norm of this piece can be made arbitrarily small by taking $r$ large (once $D$ has been shrunk to deal with the $C^0$ term).

Finally, since $A$ and $B$ are independent of $r$, the $C^k(g'_r)$-norms of the tensors they represent are bounded as $r\to\infty$. So $r^{-1/2}A$, $r^{-1}B \to 0$ in $C^k(g'_r)$, which proves the theorem.
\end{proof}

\subsection{Analysis in the local model}

This section states various analytic results concerning the K\"ahler product $S \times \C$. Since this manifold is not compact the results are not standard \emph{per se}. The situation is almost identical, however, to that of ``tubes'' considered in instanton Floer homology. A tube is the four-manifold $Y \times \R$, where $Y$ is a compact three-manifold. The proofs of the following Sobolev inequalities (as given in Chapter 3 of \cite{donaldson:fhgiymt}) carry over almost verbatim.

\begin{lemma}\label{uniform Sobolev inequalities on a product}
For indices $k$,$l$, and $q \geq p$ satisfying $k-4/p \geq l - 4/q$ there is a constant $c$ (depending only on $p$, $q$, $k$ and $l$) such that for all $\phi \in L^p_k(S \times \C)$,
$$
\| \phi \|_{L^q_l} \leq c \| \phi \|_{L^p_k}.
$$

For indices $p$, $k$ satisfying $k-4/p >0$, there is a constant $c$ depending only on $k$ and $p$, such that for all $\phi \in L^p_k(S \times \C)$,
$$
\|\phi\|_{C^0} \leq c \| \phi \|_{L^p_k}.
$$
\end{lemma}

The K\"ahler structure on $S \times \C$ determines a scalar curvature map on K\"ahler potentials. Denote the linearisation of this map by $L'$. Again, using the same arguments as in Chapter 3 of \cite{donaldson:fhgiymt} gives the following elliptic estimate for $L'$.
(The elliptic operator considered in \cite{donaldson:fhgiymt} is not $L'$; the arguments given there apply, however, to any elliptic operator determined by the local geometry.)

\begin{lemma}\label{uniform elliptic estimate for L'}
There exists a constant $A$ such that for all $\phi \in L^2_{k+4}(S \times \C)$,
$$
\| \phi\|_{L^2_{k+4}}
\leq 
A \left(\|\phi\|_{L^2} + \|L'(\phi)\|_{L^2_k}\right).
$$
\end{lemma}

\begin{lemma}\label{L'_r(uphi) - uL'_r(phi)}
There exists a constant $C$ such that for any compactly supported $u \in C^{k+4}(\C)$, and any $\phi \in L^2_{k+4}(S \times \C)$,
\begin{eqnarray*}
\|L'(u\phi) - uL'(\phi)\|_{L^p_k}
&\leq & 
C\sum_{j=1}^{k+4}\|\nabla^j u \|_{C^0}\,
\|\phi\|_{L^p_{k+4}},\\
\end{eqnarray*}
\end{lemma}

\begin{proof}
This follows from the fact that the coefficients of $L'$ are constant in the $\C$ directions.
\end{proof}

\subsection{Local analysis for $\omega_{r,n}$}

This section explains how to use Theorem \ref{local model theorem} to convert the results over $S \times \C$ from above to uniform estimates over $(X, J, \omega_{r,n})$.

First notice that, by Theorem \ref{local model theorem} with $\varepsilon <1$, over $X|_{D}$ $g_{r}-g'_{r}$ is uniformly bounded in $C^{k}(g'_{r})$. Moreover the choice of $\varepsilon$ ensures the metrics are sufficiently close that the difference $g^{-1} - g'^{-1}$ in induced metrics on the cotangent bundle is also uniformly bounded. This means that the Banach space norms on tensors determined by $g_{r}$ and $g'_{r}$ are uniformly equivalent. An immediate application of this is the following.

\begin{lemma}\label{behaviour of g_r,n-norms}
For a tensor $\alpha \in C^k(T^{*}X^{\otimes i})$, $\|\alpha\|_{C^k(g_r)} = O(1)$. Moreover, if $\alpha$ is pulled up from the base, $\|\alpha\|_{C^k(g_r)} = O\left(r^{-i/2}\right)$.
\end{lemma}

\begin{proof}
By Lemma \ref{behaviour of C^k(g'_r)}, the result is true for  the local model. Let $D$ be a disc over which Theorem \ref{local model theorem} applies for, say, $\varepsilon = 1/2$. Since $C^k(g_{r})$ and $C^k(g'_r)$ are uniformly equivalent over $X|_D$, the result holds for $C^k(g_r)$ over $X|_D$. Cover $\Sigma$ with finitely many discs $D_i$. The result holds for $C^k(g_r)$ over each $X|_{D_i}$ and hence over all of $X$.
\end{proof}

\begin{lemma}\label{approximate solution in Banach space}
~
\vspace{-1.5\baselineskip}
$$
\Scal(\omega_{r,n}) 
= 
O\left(r^{-n-1}\right) 
\quad 
\text{in $C^k(g_r)$ as $r \to \infty$,}
$$
$$
\Scal(\omega_{r,n})
= 
O\left(r^{-n-1/2}\right) \quad \text{in $L^2_k(g_r)$ as $r \to \infty$}.
$$
\end{lemma}

\begin{proof}
The expansions in negative powers of $r$ in Chapter \ref{approximate solutions} all arise via absolutely convergent power series and algebraic manipulation. This means that with respect to a \emph{fixed} metric $g$, 
$$
\Scal(\omega_{r,n}) 
= 
O\left(r^{-n-1}\right) \quad \text{in\ }C^k(g)\ \text{as\ }r\to\infty.
$$
For example, $\log(1 + r^{-1}\theta)$ is $O(r^{-1})$ in $C^k(g)$ because
\begin{eqnarray*}
\| \log(1 + r^{-1}\theta)\|_{C^k}
&\leq &
\sum_{j \geq1} r^{-(j+1)}C^j\frac{\|\theta\|^j_{C^k}}{j},\\
& = &
\log \left(1 +Cr^{-1}\|\theta\|_{C^k}\right).
\end{eqnarray*}
where $C$ a constant such that $\| \phi \psi\|_{C^k} \leq C\|\phi\|_{C^k}\|\psi\|_{C^k}$.  

The same is true with respect to the $C^k(g_r)$-norm provided that the $C^k(g_r)$-norm of a fixed function is bounded as $r \to \infty$. (Notice that the constant $C$ above does not depend on $g$.) The $C^{k}$ result now follows from Lemma \ref{behaviour of g_r,n-norms}.

To deduce the result concerning $L^2_k$-norms, notice that the $g'_r$-volume form is $r$ times a fixed form. Hence, over a disc $D$ where Theorem \ref{local model theorem} applies with $\varepsilon = 1/2$, the $g_r$-volume form is $O(r)$ times a fixed form. So the volume of $X|_D$ with respect to $g_r$ is $O(r)$. Cover $\Sigma$ with finitely many such discs, $D_i$. The volume $\vol_r$ of $X$, with respect to $g_r$, satisfies $\vol_r \leq \sum\vol(X|_{D_i}) = O(r)$. The result follows from this, the $C^{k}$ result and the fact that $\|\phi\|_{L^2_k(g_{r})}\leq (\vol_{r})^{1/2}\|\phi\|_{C^k(g_{r})}$.
\end{proof}

To transfer other results from the product to $X$, a slightly more delicate patching argument is required. Fix $\varepsilon < 1$ and cover $\Sigma$ in discs $D_1, \ldots, D_N$ satisfying the conclusions of Theorem \ref{local model theorem}. Let $\chi_i$ be a partition of unity subordinate to the cover $D_i$.

Let $\phi \in L^p_k$. Then, by the Leibniz rule and the boundedness of $\|\chi_{i}\|_{L^{p}_{k}(g_{r})}$, there exists a constant $a$ such that, for any $i = 1, \ldots , N$,
\begin{equation}\label{chi is controlled}
\| \chi_i \phi \|_{L^p_k(g_{r})} 
\leq 
a \|\phi\|_{L^p_k(g_{r})}.
\end{equation}

Everything is now in place to transfer estimates from $S \times \C$ to $(X,J,\omega_{r,n})$.

\begin{lemma}\label{uniform 1}
For indices $k$, $l$, and $q \geq p$ satisfying $k-4/p \geq l-4/q$ there is a constant $c$ (depending only on $p$, $q$, $k$ and $l$) such that for all $\phi \in L^p_k$ and all sufficiently large $r$,
$$
\| \phi \|_{L^q_l(g_{r})} \leq c \| \phi \|_{L^p_k(g_{r})}.
$$

For indices $p$, $k$ satisfying $k - 4/p\geq0$ there is a constant $c$ (depending only on $k$ and $p$) such that for all $\phi \in L^p_k$ and all sufficiently large $r$,
$$
\| \phi\|_{C^0} \leq c \| \phi \|_{L^p_k(g_{r})}.
$$
\end{lemma}

\begin{proof}
Recall the analogous result for $S \times \C$ (Lemma \ref{uniform Sobolev inequalities on a product}). Using the partition of unity $\chi_i$ from above,
$$\| \phi \|_{L^q_l(g_{r})}
\leq
\sum \| \chi_i \phi \|_{L^q_l(g_{r})}
\leq
\text{const.} \sum \| \chi_i \phi \|_{L^q_l(g'_r)},
$$
(using the uniform equivalence of the $g_{r}$- and $g'_r$-Sobolev norms). 

Considering $\chi_i \phi$ as a function over $S \times \C$, Lemma \ref{uniform Sobolev inequalities on a product} gives
$$
\| \chi_i \phi \|_{L^q_l(g'_r)}
\leq 
\text{const.}\| \chi_i \phi \|_{L^p_k(g'_r)}.
$$
Using the uniform equivalence of the  $g_{r}$- and $g'_r$-Sobolev norms again gives
$$
\|\chi_i \phi \|_{L^p_k(g'_r)}
\leq
\text{const.}\|\chi_i \phi \|_{L^p_k(g_{r})}.
$$
Finally, combining these inequalities and inequality (\ref{chi is controlled}) gives
$$
\| \phi \|_{L^q_l(g_{r})} 
\leq 
\text{const.} \sum \|\chi_i \phi \|_{L^p_k(g_{r})}
\leq 
\text{const.}\|\phi \|_{L^p_k(g_{r})}.
$$

The second Sobolev inequality is proved similarly.
\end{proof}

\begin{lemma}\label{uniform elliptic estimate for L_r,n}
There is a constant $A$, depending only on $k$, such that for all $\phi \in L^2_{k+4}$ and all sufficiently large $r$,
$$
\|\phi\|_{L^2_{k+4}(g_{r})} 
\leq 
A\left(\|\phi\|_{L^2(g_{r})} + \| L_{r}(\phi)\|_{L^2_k(g_{r})} \right).
$$
\end{lemma}

\begin{proof}
Recall the analogous result for $L'$ over $S \times \C$ (Lemma \ref{uniform elliptic estimate for L'}). This time the patching argument must be combined with Lemma \ref{Ck continuity of derivative} on the uniform continuity of the linearisation of the scalar curvature map. To apply this result, it is necessary to observe that the curvature tensor of $g'_r$ is bounded in $C^k(g'_r)$. Also, $\varepsilon$ must be taken suitably small in Theorem \ref{local model theorem}.

Using the uniform equivalence of $g'_r$- and $g_{r}$-Sobolev norms,
\begin{eqnarray*}
\| \phi \|_{L^2_{k+4}(g_{r})}
&\leq&\text{const.}
\sum \| \chi_i \phi \|_{L^2_{k+4}(g'_r)},\\
&\leq&\text{const.}
\sum \left( \| \phi \|_{L^2(g'_r)} + \|L'(\chi_i\phi)\|_{L^2_k(g'_r)}\right).
\end{eqnarray*}

Since the $\chi_i$ are functions on the base, by Lemmas  \ref{behaviour of C^k(g'_r)} and \ref{L'_r(uphi) - uL'_r(phi)},
$$
\|L'(\chi_i\phi) - \chi_i L'(\phi)\|_{L^2_k(g'_r)}
\leq
\text{const.}r^{-1/2}\|\phi\|_{L^2_k(g'_r)}.
$$
Using this, the uniform equivalence of $g'_r$- and $g_{r}$-Sobolev norms, and Lemma \ref{Ck continuity of derivative} to replace $L'$ with $L_r$ gives
$$
\| \phi \|_{L^2_{k+4}(g_{r})}
\leq
\text{const.}\left( 
\|\phi\|_{L^2(g_{r})} 
+ \| \phi \|_{L^2_k(g'_r)} + \|L_{r}(\phi)\|_{L^2_k(g_{r})} \right).
$$
This proves the result for $k=0$. It also provides the inductive step giving the result for all $k$.
\end{proof} 

\section{Global analysis}\label{inverse}

Recall that $L^2_{k,0}$ is the Sobolev space of functions with $g_{r}$-mean value zero, whilst $p$ is projection onto such functions. The aim of this section is to prove the following result.

\begin{theorem}\label{bounded inverse}
For all large $r$ and $n\geq 3$, the operator $pL_{r} : L^2_{k+4,0}\to L^2_{k,0}$ is a Banach space isomorphism. There exists a constant $C$, such that for all large $r$ and all $\psi \in L^2_{k,0}$, the inverse operator $P_{r}$ satisfies
$$
\| P_{r}\psi \|_{L^2_{k+4}(g_{r,n})} 
\leq 
Cr^{3}\| \psi \|_{L^2_k(g_{r,n})}.
$$
\end{theorem}

Unlike the uniform local results of the previous section, controlling the inverse $P_{r}$ is a global issue; indeed it is only because of global considerations (compactness of $X$, no holomorphic vector fields) that such an inverse exists. This means that the local model used in the previous chapter is not directly useful. Instead a global model is used to make calculations more straightforward.

\subsection{The global model}

Define a Riemannian metric $h_r$ on $X$ by using the fibrewise metrics determined by $\omega_0$ on the vertical vectors, and the metric $r\omega_\Sigma$ on the horizontal vectors. The metric $h_{r}$ is a Riemannian submersion on $X \to (\Sigma, r\omega_{\Sigma})$.

By construction, $g_{r,0} = h_r + a $ for some purely horizontal tensor $a \in s^{2}(T^*X)$, independent of $r$ (which is essentially given by the horizontal components of $\omega_0$). Since horizontal $1$-forms scale by $r^{-1/2}$ in the metric $h_r$ it follows immediately that for all $r$ sufficiently large,
\begin{equation}\label{global model bound}
\| g_{r,0} - h_r \|_{C^0(h_r)} \leq 1/2.
\end{equation}
Moreover, since 
$
\| g_{r} - g_{r,0} \|_{C^0(h_r)} = O\left(r^{-1}\right)
$
, inequality (\ref{global model bound}) holds with $g_{r,0}$ replaced by $g_{r}$. In particular this means that the difference in the induced metrics on the cotangent bundle is uniformly bounded and so the $L^2$-norms on tensors determined by $h_r$ and $g_{r}$ are uniformly equivalent:

\begin{lemma}\label{norm equivalence}
Let $E$ denote any bundle of tensors. There exist positive constants $k$ and $K$ such that for all $t \in \Gamma(E)$ and all sufficiently large $r$,
$$
k \| t \| _{L^2(h_r)} \leq \| t \|_{L^2(g_{r})} \leq K \| t \|_{L^2(h_r)}.
$$
\end{lemma}

\subsection{The lowest eigenvalue of $\mathscr D^* \mathscr D$}

It is more convenient to work first with the positive self-adjoint elliptic operator $\mathscr D^* \mathscr D$. Here $\mathscr D = \delb \circ \nabla$ where $\delb$ is the $\delb$-operator of the holomorphic tangent bundle and $\mathscr D^*$ is the $L^2$-adjoint of $\mathscr D$. Recall equation (\ref{derivative of scalar curvature 2}) which relates $\mathscr D^{*}\mathscr D$ to $L$.

Notice that $\mathscr D^* \mathscr D$ depends on $\omega_{r,n}$, and so on $r$. This section finds a lower bound for its first non-zero eigenvalue.

\begin{lemma}\label{no vector fields}
There are no nonzero holomorphic vector fields on $X$.
\end{lemma}

\begin{proof}
The fibres and base of $X$ have high genus. The short exact sequence of holomorphic bundles
$$
0 \to V \to TX \to \pi^*T\Sigma \to 0
$$
gives a long exact sequence in cohomology
$$
0 \to H^0(X, V) \to H^0(X, TX) \to H^0(X, \pi^*T\Sigma)\to \cdots
$$

$H^0(X, V)=0$ as the fibres admit no nonzero holomorphic vector fields. Similarly, $H^{0}(X, \pi^{*}T\Sigma) = H^{0}(\Sigma, \pi_{*}\pi^{*}T\Sigma) = H^{0}(\Sigma, T\Sigma) = 0$. The result now follows from the long exact sequence.
\end{proof}

\begin{corollary}\label{ker D = R}
$\ker \mathscr D^*\mathscr D = \R$. Equivalently, $\mathscr D^* \mathscr D \colon L^2_{k+4,0} \to L^2_{k,0}$ is an isomorphism.
\end{corollary}

\begin{proof}
$\ker \mathscr D^*\mathscr D = \ker \mathscr D$ is those functions with holomorphic gradient. The previous lemma implies such functions must be constant. The second statement now follows from the fact that $\mathscr D^* \mathscr D$ is a self-adjoint index zero operator.
\end{proof}

To find a lower bound for the first non-zero eigenvalue of $\mathscr D^* \mathscr D$, similar bounds are first found for the Hodge Laplacian and for the $\delb$-Laplacian on sections of the holomorphic tangent bundle.

\begin{lemma}\label{first eval of laplacian}
There exists a positive constant $C_1$ such that for all $\phi$ with $g_{r}$-mean value zero and all sufficiently large $r$, 
$$
\left\| \diff \phi \right\|^2_{L^2(g_{r})} 
\geq 
C_1r^{-1} \left\| \phi \right\|^2_{L^2(g_{r})}.
$$
\end{lemma}

\begin{proof}
There exists a constant $m$ such that $\phi - m$ has $h_1$-mean value zero. Since $m$ is constant, $\diff\phi = \diff(\phi -m)$. Using Lemma \ref{norm equivalence}, 
$
\left\| \diff \phi \right\|_{L^2(g_{r})} 
\geq \text{const.}
\left\| \diff (\phi-m) \right\|_{L^2(h_r)}
$.

Let $|\cdot|_{h_r}$ denote the pointwise inner product defined by $h_r$. By definition of $h_r$ it follows that 
$
\left| \diff(\phi - m) \right|^2_{h_r} 
\geq
r^{-1}\left| \diff( \phi -m) \right|^2_{h_1}
$.
Moreover, the volume forms satisfy $\diff\!\vol(h_r) = r \diff\!\vol(h_1)$. Hence,
$$
\left\| \diff (\phi - m) \right\|^2_{L^2(h_r)} 
\geq 
\left \|\diff(\phi -m)\right\|^2_{L^2(h_1)}.
$$

Now $\phi-m$ has $h_1$-mean value zero. Let $c$ be the first eigenvalue of the $h_1$-Laplacian. Then
$$\left \|\diff(\phi -m)\right\|^2_{L^2(h_1)} 
\geq 
c \left\|\phi -m\right\|^2_{L^2(h_1)}
 =  cr^{-1} \left\| \phi - m \right\|^2_{L^2(h_r)}.
$$
Using Lemma \ref{norm equivalence} again gives
$$
\left \| \phi - m \right\|^2_{L^2(h_r)}
\geq
\text{const.} \left \| \phi - m \right\|^2_{L^2(g_{r})}\\
\geq
\text{const.} \left\| \phi  \right\|^2_{L^2(g_{r})}
$$
where the second inequality follows from the fact that $\phi$ has $g_{r}$-mean value zero. 

Putting the pieces together completes the proof.
\end{proof}

\begin{lemma}\label{first eval of delb laplacian}
There exists a positive constant $C_2$ such that for all $\xi \in \Gamma(TX)$ and all sufficiently large $r$, 
$$
\left\| \delb \xi \right\|^2_{L^2(g_{r})} 
\geq 
C_2r^{-2} \left\| \xi \right\|^2_{L^2(g_{r})}.
$$
\end{lemma}

\begin{proof}
The proof is similar to that of Lemma \ref{first eval of laplacian} above. By Lemma \ref{norm equivalence},
$$
\left\| \delb \xi \right\|^2_{L^2(g_{r})} 
\geq 
\text{const.}\left\| \delb \xi \right\|^2_{L^2(h_r)}.
$$
By definition of $h_{r}$,
$
\left| \delb \xi \right|^2_{h_r} \geq r^{-1} \left| \delb \xi \right|^2_{h_1}
$. 
Using this and $\diff\!\vol(h_r) = r\diff\!\vol(h_1)$ gives
$
\left\| \delb \xi \right\|^2_{L^2(h_r)} 
\geq 
\left\| \delb \xi \right\|^2_{L^2(h_1)}
$.

Let $c$ be the first eigenvalue of the $\delb$-Laplacian determined by the metric $h_1$. Then
$
\left\| \delb \xi \right\|^2_{L^2(h_1)}
\geq
c \left\| \xi \right\|^2_{L^2(h_1)}
$.
By definition of $h_{r}$,
$
\left| \xi \right|^2_{h_1} \geq r^{-1} \left| \xi \right|^2_{h_r}
$.
Hence,
$
\left\| \xi \right\|^2_{L^2(h_1)} \geq r^{-2}\left\| \xi \right\|^2_{L^2(h_r)}
$.
Finally, using Lemma \ref{norm equivalence} to convert back to the $L^2(g_{r})$-norm of $\xi$, and putting all the pieces together gives the result.
\end{proof}

\begin{lemma}\label{first eval of D^*D}
There exists a constant $C$ such that for all $\phi$ with $g_{r}$-mean value zero and all sufficiently large $r$,
$$
\left\| \mathscr D \phi \right\|^2_{L^2(g_{r})} 
\geq 
Cr^{-3} \left\| \phi \right\|^2_{L^2(g_{r})}.
$$
\end{lemma}

\begin{proof}
Combining Lemmas \ref{first eval of laplacian} and \ref{first eval of delb laplacian} shows that whenever $\phi$ has $g_{r}$-mean value zero,
\begin{eqnarray*}
\left\| \delb \nabla \phi \right\|^2_{L^2(g_{r})}
&\geq&
C_2r^{-2} \| \nabla \phi \|^2_{L^2(g_{r})},\\
&=&
C_2r^{-2} \| \diff \phi \|^2_{L^2(g_{r})},\\
&\geq&
C_1C_2r^{-3} \| \phi \|^2_{L^2(g_{r})}.
\end{eqnarray*}
\end{proof}

\subsection{A uniformly controlled inverse}
 
\begin{lemma}\label{elliptic estimate for D^*D}
There is a constant $A$, depending only on $k$, such that for all $\phi \in L^2_{k+4}$ and sufficiently large $r$,
$$
\|\phi\|_{L^2_{k+4}(g_{r})} 
\leq 
A\left(\|\phi\|_{L^2(g_{r})} 
+ 
\| \mathscr D^* \mathscr D(\phi)\|_{L^2_k(g_{r})} \right).
$$
\end{lemma}

\begin{proof}
Recall equation (\ref{derivative of scalar curvature 2}): $
L_{r}(\phi) = \mathscr D^* \mathscr D (\phi) + 
\nabla \Scal(\omega_{r,n}) \cdot \nabla \phi
$. Since $\Scal(\omega_{r,n})$ tends to zero in $C^k(g_{r})$, $L_{r} - \mathscr D^* \mathscr D$ converges to zero in operator norm calculated with respect to the $L^2_k(g_{r})$-Sobolev norms. Hence the estimate follows from the analogous result for $L_{r}$ (Lemma \ref{uniform elliptic estimate for L_r,n}).
\end{proof}

\begin{theorem}\label{bounded inverse for D^*D}
The operator
$
\mathscr D^* \mathscr D : L^2_{k+4,0}\to L^2_{k,0}
$
is a Banach space isomorphism. There exists a constant $K$, such that for all large $r$ and all $\psi \in L^2_{k,0}$, the inverse operator $Q_r$ satisfies
$$
\| Q_r\psi \|_{L^2_{k+4}(g_{r})} 
\leq 
Kr^{3}\| \psi \|_{L^2_k(g_{r})}.
$$
\end{theorem}

\begin{proof}
The inverse $Q_r$ exists by Corollary \ref{ker D = R}. It follows from Lemma \ref{first eval of D^*D} applied to $\phi = Q_r \psi$ that there is a constant $C$ such that for all $\psi \in L^2_{k,0}$,
$$
\| Q_r\psi \|_{L^2(g_{r})} \leq Cr^3 \| \psi\|_{L^2(g_{r})}.
$$
Applying Lemma \ref{elliptic estimate for D^*D} to $\phi = Q_r\psi$ extends this bound to the one required.
\end{proof}

Next recall the following standard result (proved via a geometric series).

\begin{lemma}\label{invertibility is open}
Let $D \colon B_1 \to B_2$ be an bounded invertible linear map of Banach spaces with bounded inverse $Q$. If $L \colon B_1 \to B_2$ is another linear map with 
$$
\| L - D \| \leq (2\|Q\|)^{-1},
$$ 
then $L$ is also invertible with bounded inverse $P$ satisfying 
$
\| P \| \leq 2\|Q\|
$.
\end{lemma}

The pieces are now in place to prove Theorem \ref{bounded inverse} (which is stated at the start of this chapter).
 
\begin{proof}[Proof of Theorem \ref{bounded inverse}]
Since 
$(L_{r} - \mathscr D^* \mathscr D)\phi 
= 
\nabla \Scal(\omega_{r,n}) \cdot \nabla \phi$,
there exists a constant $c$ such that, in operator norm computed with respect to the $g_{r}$-Sobolev norms,
$\| pL_{r} - \mathscr D^* \mathscr D \| \leq cr^{-n-1}$.

So for $n\geq3$, and for large enough $r$, \
$\| pL_{r} - \mathscr D^* \mathscr D \| \leq (2\|Q_r\|)^{-1}$.
Lemma \ref{invertibility is open} shows that $pL_{r}$ is invertible and gives the upper bound 
$$
\| P_{r} \| \leq 2\|Q_r\| \leq Cr^3
$$
for some $C$.
\end{proof}

\subsection{An improved bound}

It should be possible to improve on this estimate. For example, over a  product of two curves a simple separation of variables argument shows that $\|P_r\| = Cr^2$.

The above proof of Theorem \ref{bounded inverse} concatenates two eigenvalue estimates, each of which is saturated only when applied to an eigenvector corresponding to the first eigenvalue. Certainly over a product, the functions which get closest to saturating the first estimate (Lemma \ref{first eval of laplacian}) have gradients which can be controlled more efficiently than is done in the proof of the second estimate (Lemma \ref{first eval of delb laplacian}). 

In general, it should be possible to obtain a better bound for $\| P_{r}\|$ by examining this interplay between the two eigenvalue estimates. However, the bound proved above is sufficient to complete the proof of Theorem \ref{main theorem}. 

\section{Loose ends}\label{loose ends}

\subsection{Controlling the nonlinear terms}

Denote by $\Scal_{r}$ the scalar curvature map on K\"ahler potentials determined by $\omega_{r,n}$: $\Scal_{r}(\phi) = \Scal(\omega_{r} + i\delb\del \phi)$. Recall that $S_{r} = p\Scal_{r}$. Denote by $N_{r} = S_{r} - pL_{r}$ the nonlinear terms of $S_{r}$

\begin{lemma}\label{controlling N_r,n}
Let $k\geq3$. There exists positive constants $c$ and $K$, such that for all $\phi$, $\psi \in L^2_{k+4}$ with $\|\phi\|_{L^2_{k+4}}$, $\|\psi\|_{L^2_{k+4}} \leq c$ and for sufficiently large $r$,
$$
\|N_{r}(\phi) - N_{r}(\psi)\|_{L^2_k}
\leq
K \max\left\{ \|\phi\|_{L^2_{k+4}},\,\|\psi\|_{L^2_{k+4}} \right\}
\| \phi - \psi \|_{L^2_{k+4}}$$
where $g_{r}$-Sobolev norms are used throughout.
\end{lemma}

\begin{proof}
By the mean value theorem, 
$$
\| N_{r}(\phi) - N_{r}(\psi) \|_{L^2_k(g_{r})} 
\leq 
\sup_{\chi \in [\phi,\psi]} 
\|(DN_{r})_\chi\| \| \phi - \psi \|_{L^2_{k+4}(g_{r})} 
$$
where $(DN_{r})_\chi$ is the derivative of $N_{r}$ at $\chi$. 

Now $DN_{r} = p(L_{r})_\chi - p L_{r}$ where $(L_{r})_\chi$ is the linearisation of $\Scal_{r}$ at $\chi$. In other words, $(L_{r})_{\chi}$ is the linearisation of the scalar curvature map determined by the metric $\omega_{r,n}+i\delb\del\chi$. Applying Lemma \ref{Lpk continuity of derivative} to this metric and $\omega_{r,n}$ gives 
$
\|(L_{r})_\chi - L_{r} \| 
\leq 
\text{const.} \|\chi\|_{L^2_{k+4(}g_{r})}
$. 
As $k\geq3$ the condition on the indices in Lemma \ref{Lpk continuity of derivative} is met. Notice also that Lemma \ref{Lpk continuity of derivative} requires the constants in the $g_{r}$-Sobolev inequalities to be uniformly bounded --- which is proved in Lemma \ref{uniform 1} --- and the $C^k(g_{r})$-norm of the curvature of $\omega_{r,n}$ to be bounded above --- which follows from Theorem \ref{local model theorem} and Lemma \ref{Ck continuity of curvature}.

Since $p$ is uniformly bounded (an $L^2_k(g_{r})$-orthogonal projection even) and since, for all $\chi \in [\phi, \psi]$, 
$
\| \chi \|_{L^{2}_{k+4}} 
\leq 
\max\{ \|\phi\|_{L^{2}_{k+4}},\,\|\psi\|_{L^{2}_{k+4}} \}
$
the result follows.
\end{proof} 

\subsection{Completing the proof}

\begin{proof}[Proof of Theorem \ref{main theorem}]
For all large $r$ and $n \geq 3$, the map 
$$
S_{r}
\colon 
L^2_{k+4,0}(g_{r})
\to 
L^2_{k,0}(g_{r})
$$ 
has the following properties:

\begin{enumerate}
\item
$S_{n}(0) = O\left(r^{-n -1/2}\right)$ in $L^2_k(g_{r})$, by Lemma \ref{approximate solution in Banach space} and the fact that $p$ has operator norm $1$.
\item
The derivative of $S_{r}$ at the origin is an isomorphism with inverse $P_{r}$ which is $O(r^3)$. This is proved in Theorem \ref{bounded inverse}.
\item
There exists a constant $K$ such that for all sufficiently small $M$, the nonlinear piece $N_{r}$ of $S_{r}$ is Lipschitz with constant $M$ on a ball of radius $KM$. This follows directly from Lemma \ref{controlling N_r,n}.
\end{enumerate}

Recall the statement of the inverse function theorem (Theorem \ref{ift}). The second and third of the above properties imply that the radius $\delta'_{r}$ of the ball about the origin on which $N_{r}$ is Lipschitz with constant $(2\|P_{r}\|)^{-1}$ is bounded below by $Cr^{-3}$ for some positive $C$. As $\delta_{r} = \delta'_{r}(2\|P_{r}\|)^{-1}$, it follows that $\delta_{r}$ is bounded below by $Cr^{-6}$ for some positive $C$.

Hence for $\psi \in L^{2}_{k}$ with $\|S_{r}(0) - \psi \|_{L^{2}_{k}(g_{r})}\leq Cr^{-6}$ the equation $S_{r}(\phi) = \psi$ has a solution. In particular, the first of the above properties implies that, for $n\geq6$ and sufficiently large $r$, the equation $S_{r}(\phi) = 0$ has a solution.

Since $\Scal$ differs from $S_{r}$ by a constant, the metric $\omega_{r,n}+i\delb\del \phi$ has constant scalar curvature. Iteratively applying the regularity Lemma \ref{regularity} (which can be done provided $k$ is high enough to ensure that $L^{2}_{k+4} \hookrightarrow C^{2,\alpha}$) gives that $\phi$ is smooth.
\end{proof}

\section{Higher dimensional varieties}\label{higher dimensions}

Theorem \ref{main theorem} extends to certain higher dimensional fibrations. The required conditions are set out below. To understand their relevance, compare  the summary at the end of Chapter \ref{approximate solutions}.

\begin{itemize}
\item[(A)]
\emph
{Let X be a compact connected K\"ahler manifold with no nonzero holomorphic vector fields and $\pi\colon X\to B$ a holomorphic submersion.}
\end{itemize}

Let $\kappa_{0}$ be a K\"ahler class on $X$; denote by $\kappa_{b}$ the K\"ahler class on the fibre $F_{b}$ over $b$ obtained by restricting $\kappa_{0}$.

\begin{itemize}
\item[(B)]
\emph
{For every $b\in B$, $\kappa_{b}$ contains a unique cscK metric $\omega_{b}$; the form $\omega_{b}$ depends smoothly on $b$.}
\end{itemize}

Let $\omega$ be a K\"ahler form representing $\kappa_{0}$. For each $b$ there is a unique function $\phi_{b}\in C^{\infty}(F_{b})$ with $\omega_{b}$-mean value zero such that the fibrewise restriction of $\omega$ plus $i\delb\del\phi_{b}$ is $\omega_{b}$. The smoothness assumption in (B) implies that the $\phi_{b}$ fit together to give a smooth function $\phi \in C^{\infty}(X)$; so $\omega_{0}=\omega+i\delb\del\phi$ is a K\"ahler metric in $\kappa_{0}$ whose fibrewise restriction has constant scalar curvature.

The metric $\omega_{0}$ gives a Hermitian structure in the vertical tangent bundle $V$ and hence also in the line bundle $\Lambda^{\text {max}}V^{*}$; denote its curvature $F$. Taking the fibrewise mean value of the horizontal-horizontal component of $iF$ (with respect to the metric $\omega_{0}$) defines a form $\alpha \in \Omega^{1,1}(B)$.

\begin{itemize}
\item[(C)]
\emph
{There is a metric $\omega_{B}$ on the base solving $\Scal(\omega_{B}) - \Lambda \alpha = {\text const}$; there are no nontrivial deformations of $\omega_{B}$ through cohomologous solutions to this equation.}
\end{itemize}

\begin{theorem}\label{higher dimensional theorem}
Let $X$ satisfy (A), (B) and (C). Then for all large $r$ the K\"ahler class 
$
\kappa_r = \kappa_0 + r[\omega_B]
$
contains a constant scalar curvature K\"ahler metric.
\end{theorem}

\begin{proof}
The proof follows the same lines as that of Theorem \ref{main theorem}. Using the same notation as in the preceding paragraphs, let $\omega_{r} = \omega_{0} + r\omega_{B}$, where $\omega_{B}$ is the solution to the PDE mentioned in (C). Notice that $[\omega_{r}] = \kappa_{r}$.

The first step is to construct approximate solutions of arbitrary accuracy. An identical calculation to that in Lemma \ref{O(r^-1) expansion of Scal(omega_r)} shows that $\Scal(\omega_{r}) = \Scal(\omega_{b}) + \sum r^{-j}\psi_{j}$ for some functions $\psi_{j}$; moreover the fibrewise mean value of of $\psi_{1}$ is $\Scal(\omega_{B}) - \Lambda_{\omega_{B}}\alpha$. By assumptions (B) and (C) then, the fibrewise mean value of $\Scal(\omega_{r})$ is constant to $O(r^{-2})$. 

The same argument which proves Lemma \ref{L_r to O(r^-1)} gives that the linearisation of the scalar curvature map determined by $\omega_{r}$ satisfies $L_{r} = L_{b} + O(r^{-1})$, where $L_{b}$ is the linearisation of the scalar curvature map on the fibre over $b$ determined by $\omega_{b}$. By assumption (B), $\omega_{b}$ is the unique constant scalar curvature metric in $\kappa_{b}$, hence $\ker L_{b}$ is the constant functions on $F_{b}$. This, elliptic regularity and the fact that $L_{b}$ is self adjoint mean that for any $\Theta \in C^{\infty}(F_{b})$ with $\omega_{b}$-mean value zero, there exists a unique $\phi \in C^{\infty}(F_{b})$ with $\omega_{b}$-mean value zero satisfying $L_{b}\phi = \Theta$. 

Applying this argument fibrewise, as in the proof of Lemma \ref{invertibility of O(r^-1) piece of L_r}, shows that $L_{b}$ is a bijection on $C^{\infty}_{0}(X)$ (\emph{i.e.}\ on functions with fibrewise mean value zero). This guarantees the existence of a function $\phi_{1} \in C^{\infty}_{0}(X)$ with 
$$
\Scal\left(\omega_{r}+i\delb\del r^{-1}\phi_{1}\right) 
= 
c_{0} +c_{1}r^{-1} + \sum_{j \geq 2}r^{-j}\chi_{j}
$$
for some constants $c_{0}$ and $c_{1}$ and functions $\chi_{j}$. 

The higher order approximate solutions are constructed exactly as in Sections \ref{third} and \ref{higher}. The argument hinges on the surjectivity of two particular linear differential operators. The first is the operator $L_{b}\colon C^{\infty}_{0}(X) \to C^{\infty}_{0}(X)$ whose surjectivity is justified above. The second is the linearisation of the map $C^{\infty}(B) \to C^{\infty}(B)$ given by 
$$
G \colon f 
\mapsto 
\Scal(\omega_{B} + i\delb\del f) 
+ 
\Lambda_{\omega_{B} + i \delb\del f}\alpha.
$$
Notice that $\int G(f) (\omega_{B}+i\delb\del f)^{\dim B}$ is independent of $f$; by assumption (C), $G(0)$ is constant; hence differentiating gives that $\im DG$ is $L^{2}$-orthogonal to the constants. By assumption (C), $\ker DG$ is the constants. This, elliptic regularity and the fact that $DG$ has index zero implies that 
$DG$ is surjective onto smooth functions on $B$ with mean value zero.

In conclusion, for any integer $n$, there exist functions $f_{1}, \ldots , f_{n-1} \in C^{\infty}(B)$, $\phi_{1}, \ldots , \phi_{n} \in C^{\infty}_{0}(X)$ and constants $c_{0}, \ldots, c_{n}$ such that the metric
$$
\omega_{r,n} 
= 
\omega_{r} 
+ 
i\delb\del \sum_{j=1}^{n-1}r^{-j}f_{j}
+
i\delb\del \sum_{j=1}^{n}r^{-j}\phi_{j}
$$
satisfies
$$
\Scal(\omega_{r,n}) = \sum_{j=0}^{n}c_{j}r^{-j} 
+ 
O\left(r^{-n-1}\right).
$$

The remainder of the proof of Theorem \ref{main theorem} is not dimension specific. Similar to Section \ref{local analysis}, the local model is given by $(F \times \C^{\dim B}, \omega_{F}\oplus\omega_{\text{flat}})$ with $\omega_{F}$ a constant scalar curvature metric on $F$. The same Sobolev inequalities (modulo the obvious changes regarding the indices) hold  over this space as over $S \times \C$. (Indeed they hold over any manifold with uniformly bounded geometry, as is remarked in \cite{donaldson:fhgiymt}.) These estimates transfer to $X$ in an identical way to before. 

The global model is, as in Section \ref{inverse}, a Riemannian submersion constructed by ignoring the horizontal contribution of $\omega_{0}$. The remaining steps in the proof now proceed identically.
\end{proof}

It is, perhaps, suprising that condition (C) is more complicated than just the existence of a constant scalar curvature metric on $B$. Indeed the equation in (C) involves the whole of $X$ and is not just a condition on the base. Bearing in mind the conjectured correspondence with stably polarised varieties, this may have an algebro-geometric interpretation. Namely, $X$ may be stable with respect to the polarisation $\kappa_0 + r[\omega_B]$ if the fibres are stably polarised by the restriction of $\kappa_0$ and if the base is also stably polarised, not with respect to $[\omega_B]$, but rather some other polarisation constructed from $[\omega_B]$ and the push down of the top exterior power of the vertical cotangent bundle of $X$. 

\bibliographystyle{plain}
\bibliography{../../bibliographies/bibthesis}

\begin{thebibliography}{10}

\bibitem{atiyah:tsofb}
M.~F. Atiyah.
\newblock The signature of fibre bundles.
\newblock In {\em Global Analysis (Papers in Honor of K. Kodaira)}, pages
  73--84. University of Tokyo Press, 1969.

\bibitem{aubin:naommae}
T.~Aubin.
\newblock {\em Nonlinear analysis on manifolds. {M}onge-{A}mp\`ere equations},
  volume 252 of {\em Grundlehren der Mathematischen Wissenschaften [Fundamental
  Principles of Mathematical Sciences]}.
\newblock Springer-Verlag, New York, 1982.

\bibitem{donaldson:rogtcga4mt}
S.~K. Donaldson.
\newblock Remarks on gauge theory, complex geometry and {$4$}-manifold
  topology.
\newblock In {\em Fields Medallists' lectures}, volume~5 of {\em World Sci.
  Ser. 20th Century Math.}, pages 384--403. World Sci. Publishing, River Edge,
  NJ, 1997.

\bibitem{donaldson:fhgiymt}
S.~K. Donaldson.
\newblock {\em Floer Homology Groups in Yang-Mills Theory}.
\newblock Cambridge Tracts in Mathematics. Cambridge University Press, 2002.

\bibitem{harris.morrison:moc}
J.~Harris and I.~Morrison.
\newblock {\em Moduli of curves}, volume 187 of {\em Graduate Texts in
  Mathematics}.
\newblock Springer-Verlag, New York, 1998.

\bibitem{hong:rmwchsc}
Y.-J. Hong.
\newblock Ruled manifolds with constant {H}ermitian scalar curvature.
\newblock {\em Math. Res. Lett.}, 5(5):657--673, 1998.

\bibitem{hong:chsceorm}
Y.-J. Hong.
\newblock Constant {H}ermitian scalar curvature equations on ruled manifolds.
\newblock {\em J. Differential Geom.}, 53(3):465--516, 1999.

\bibitem{kazdan.warner:icfpdewatrg}
J.~L. Kazdan and F.~W. Warner.
\newblock Integrability conditions for {$\Delta u=k-Ke\sp{\alpha u}$} with
  applications to {R}iemannian geometry.
\newblock {\em Bull. Amer. Math. Soc.}, 77:819--823, 1971.

\bibitem{kodaira:actoias}
K.~Kodaira.
\newblock {A certain type of irregular algebraic surface}.
\newblock {\em Journal Analyse Mathematique}, 19:207--215, 1967.

\bibitem{schechter:hoaaif}
E.~Schechter.
\newblock {\em Handbook of analysis and its foundations}.
\newblock Academic Press Inc., San Diego, CA, 1997.

\bibitem{tian:kemwpsc}
G.~Tian.
\newblock K\"ahler-{E}instein metrics with positive scalar curvature.
\newblock {\em Invent. Math.}, 130(1):1--37, 1997.

\bibitem{yau:opig}
S.-T. Yau.
\newblock Open problems in geometry.
\newblock In {\em Differential geometry: partial differential equations on
  manifolds (Los Angeles, CA, 1990)}, volume~54 of {\em Proc. Sympos. Pure
  Math.}, pages 1--28. Amer. Math. Soc., Providence, RI, 1993.

\end{thebibliography}

\end{document}